\documentclass[12pt,reqno]{amsart}
\usepackage{amsmath}
\usepackage{cite}
\usepackage[latin1]{inputenc}
\usepackage{amssymb}
\usepackage{amscd}
\usepackage{mathrsfs}
\usepackage[english]{babel}
\usepackage[babel=true]{csquotes}
\usepackage[pdftex]{graphicx}
\usepackage{graphicx}
\newcommand{\ue}[0]{\ensuremath{\mathbf{u}}}
\newcommand{\J}[0]{\ensuremath{\mathbf{J}}}
\newcommand{\ve}[0]{\ensuremath{\mathbf{v}}}
\newcommand{\ix}[0]{\ensuremath{\mathbf{x}}}

\newcommand{\ka}[0]{\ensuremath{\mathbf{k}}}
\newcommand{\en}[0]{\ensuremath{\mathbf{n}}}
\newcommand{\er}[0]{\ensuremath{\mathbf{r}}}

\newcommand{\One}[0]{\ensuremath{\mathbf{1}}}

\newcommand{\N}[0]{\ensuremath{\mathbb{N}}}

\newcommand{\E}[0]{\ensuremath{\mathbf{E}}}
\newcommand{\F}[0]{\ensuremath{\mathcal{F}}}

\newcommand{\Pe}[0]{\ensuremath{\mathbf{P}}}
\newcommand{\R}[0]{\ensuremath{\mathbb{R}}}

\newcommand{\Z}[0]{\ensuremath{\mathbb{Z}}}

\newtheorem{theo}{\sc{Theorem}}[section]
\newtheorem{lemm}[theo]{\sc{Lemma}}
\newtheorem{defi}[theo]{\sc{Definition}}

\newtheorem{cor}[theo]{\sc{Corollary}}
\newtheorem{rmq}[theo]{\sc{Remark}}

\author{Alexandre Boritchev}
\title[]{Turbulence for the generalised Burgers equation}
\date{\today}
\begin{document}


\maketitle

\textbf{Abstract.} In this survey, we review the results on turbulence for the generalised Burgers equation on the circle:
\begin{equation} \nonumber
u_t+f'(u)u_x=\nu u_{xx}+\eta,\ x \in S^1=\R/\Z,
\end{equation}
obtained by A.Biryuk and the author in \cite{Bir01,BorK,BorW,BorD}. Here, $f$ is smooth and strongly convex, whereas the constant $0<\nu \ll 1$ corresponds to a viscosity coefficient. We will consider both the case $\eta=0$ and the case when $\eta$ is a random force which is smooth in $x$ and irregular (kick or white noise) in $t$. In both cases, sharp bounds for Sobolev norms of $u$ averaged in time and in ensemble of the type $C \nu^{-\delta},\ \delta \geq 0$, with the same value of $\delta$ for upper and lower bounds, are obtained. These results yield sharp bounds for small-scale quantities characterising turbulence, confirming the physical predictions \cite{BK07}.


\section*{Abbreviations}

\begin{itemize}
\item 1d, 3d, multi-d: 1, 3, multi-dimensional
\item a.e.: almost everywhere
\item a.s.: almost surely
\item (GN): the Gagliardo--Nirenberg inequality (Lemma~\ref{GN})
\item i.i.d.: independent identically distributed
\item r.v.: random variable
\end{itemize}

\section*{Introduction} \label{intro}

The generalised 1d space-periodic Burgers equation
\begin{equation} \label{Burbegin}
\frac{\partial u}{\partial t} +  f'(u)\frac{\partial u}{\partial x}  - \nu \frac{\partial^2 u}{\partial x^2} = 0,\quad \nu>0,\ x \in S^1=\R/\Z
\end{equation}
(the classical Burgers equation \cite{Bur74} corresponds to $f(u)=u^2/2$) is a popular model for the Navier--Stokes equation. Indeed, both of them have similar nonlinearities and dissipative terms. For $\nu \ll 1$ and $f$ strongly convex, i.e. satisfying:
\begin{equation} \label{strconvex}
f''(x) \geq \sigma > 0,\quad x \in \R,
\end{equation}
solutions of (\ref{Burbegin}) exhibit turbulent-like behaviour, called \enquote{Burgulence} \cite{BF01,BK07}. To simplify the presentation, we restrict ourselves to solutions with zero mean value in space:
\begin{equation} \label{zero}
\int_{S^1}{u(t,x) dx}=0,\quad \forall t \geq 0.
\end{equation}
The space mean value is a conserved quantity. Indeed, since $u$ is $1$-periodic in space, we have
\begin{align} \nonumber
\frac{d}{dt}\int_{S^1}{u(t,x) dx}&=-\int_{S^1}{f'(u(t,x))u_x(t,x) dx}+\nu \int_{S^1}{u_{xx}(t,x) dx}=0.
\end{align}
Thus, it suffices to assume that the initial value $u_0=u(0,\cdot)$ satisfies (\ref{zero}). If the mean value of $u_0$ on $S^1$ equals $b$, we may consider the zero mean value function
$$
v(t,x)=u(t,x+bt)-b,
$$
which is a solution of (\ref{Burbegin}) with $f(y)$ replaced by $g(y)=f(y+b)-by$.
\\ \indent
In this survey, we consider both the unforced equation (\ref{Burbegin}) and the generalised Burgers equation with an additive forcing term, smooth in space and highly irregular in time (see Subsection~\ref{force}). We summarise the estimates obtained by A.Biryuk and the author \cite{Bir01,BorK,BorW,BorD} for Sobolev norms as well as for the small-scale quantities relevant for the theory of hydrodynamical turbulence (the dissipation length scale, the structure functions and the energy spectrum). This survey is partially based on the Ph.D. thesis of the author \cite{BorPhD}, where some technical points are covered in more detail.
\\ \indent
For the \textit{unforced} Burgers equation, some \textit{upper} estimates for small-scale quantities have been obtained previously. Lemma~\ref{detuxpos} is an analogue in the periodic setting of the one-sided Lipschitz estimate due to Oleinik, and the upper estimate for the structure function $S_1(\ell)$ follows from an estimate for the solution in the class of bounded variation functions $BV$. For references on these classical aspects of the theory of scalar conservation laws, and namely for Oleinik's estimate, see \cite{Daf10,Lax06,Ser99}. For some upper estimates for small-scale quantities, see \cite{Kre88,Tad93}.
\\ \indent
The research on small-scale behaviour of solutions for this nonlinear PDE is motivated by the problem of turbulence. It has been inspired by the pioneering works of Kuksin, who obtained lower and upper estimates for Sobolev norms by negative powers of the viscosity for a large class of equations (see \cite{Kuk97GAFA,Kuk99GAFA} and references in \cite{Kuk99GAFA}). For more recent results obtained by Kuksin, Shirikyan and others for the 2D Navier--Stokes equation, see the book \cite{KuSh12} and references therein. 
\\ \indent
Estimates for Sobolev norms as well as for small-scale quantities obtained here are asymptotically sharp in the sense that viscosity enters lower and upper bounds at the same negative power. Such estimates are not available for the more complicated equations considered in \cite{Kuk97GAFA,Kuk99GAFA,KuSh12}.
\\ \indent
We do not consider other aspects of Burgers turbulence, such as the inviscid limit or the behaviour of solutions for spatially rough forcing, and we refer the reader to the survey \cite{BK07}.
\medskip \\ \indent
\textbf{Organisation of the paper:} We begin by introducing the notation and setup in Section~\ref{nota}. In Section~\ref{turb}, we present the K41 theory as well as the physical predictions for Burgers turbulence. In Section~\ref{results}, we formulate the main results. 
\\ \indent
In Section~\ref{det}, we consider the solution $u(t,x)$ of the unforced equation (\ref{Burbegin}). In Subsection~\ref{detSob}, we begin by recalling the upper estimate for the quantity
$$
\max_{s \in [t,t+1],\ x \in S^1}{u_x(s,x)},\ t \geq 1.
$$
Using this bound, we get upper and lower estimates for the Sobolev norms of $u$. In Subsection~\ref{detturb} we study the implications of our results in terms of the theory of Burgulence.  Namely, we give sharp upper and lower bounds for the dissipation length scale, the increments and the spectral asymptotics for the flow $u(t,x)$. These bounds hold uniformly for $\nu \leq \nu_0$, where $\nu_0>0$ depends only on $f$ and on the initial condition. Those results rigorously justify the physical predictions for small-scale quantities.
\\ \indent
In Section~\ref{rand}, we consider the randomly forced generalised Burgers equation and we obtain analogues of the results in Section~\ref{det}, which also confirm the physical predictions \cite{BK07}. In Section~\ref{stat}, we are concerned with the stationary measure for the randomly forced generalised Burgers equation.

\section{Notation and setup} \label{nota}

All functions which we consider in this paper are real-valued.

\subsection{Functional spaces and Sobolev norms} \label{Sob}

Consider a zero mean value integrable function $v$ on $S^1$. For $p \in [1,\infty]$, we denote its $L_p$ norm by $\left|v\right|_p$. The $L_2$ norm is denoted by  $\left|v\right|$, and $\left\langle \cdot,\cdot\right\rangle$ stands for the $L_2$ scalar product. From now on $L_p,\ p \in [1,\infty],$ denotes the space of zero mean value functions in $L_p(S^1)$. Similarly, $C^{\infty}$ is the space of $C^{\infty}$-smooth zero mean value functions on $S^1$.
\\ \indent
For a nonnegative integer $m$ and $p \in [1,\infty]$, $W^{m,p}$ stands for the Sobolev space of zero mean value functions $v$ on $S^1$ with finite homogeneous norm
\begin{equation} \nonumber
\left|v\right|_{m,p}=\left|\frac{d^m v}{dx^m}\right|_p.
\end{equation}
In particular, $W^{0,p}=L_p$ for $p \in [1,\infty]$. For $p=2$, we denote $W^{m,2}$ by $H^m$ and abbreviate the corresponding norm as $\left\|v\right\|_m$. 
\\ \indent
Since the length of $S^1$ is $1$, we have
$$
|v|_1 \leq |v|_{\infty} \leq |v|_{1,1} \leq |v|_{1,\infty} \leq \dots \leq |v|_{m,1} \leq |v|_{m,\infty} \leq \dots
$$
We recall a version of the classical Gagliardo--Nirenberg inequality (see \cite[Appendix]{DG95}):
\begin{lemm} \label{GN}
For a smooth zero mean value function $v$ on $S^1$,
$$
\left|v\right|_{\beta,r} \leq C \left|v\right|^{\theta}_{m,p} \left|v\right|^{1-\theta}_{q},
$$
where $m>\beta\geq 0$, and $r$ is defined by
$$
\frac{1}{r}=\beta-\theta \Big( m-\frac{1}{p} \Big)+(1-\theta)\frac{1}{q},
$$
under the assumption $\theta=\beta/m$ if $p=1$ or $p=\infty$, and $\beta/m \leq \theta < 1$ otherwise. The constant $C$ depends on $m,p,q,\beta,\theta$.
\end{lemm}
From now on, we will refer to this inequality as (GN).
\\ \indent
For any $s \geq 0$, $H^{s}$ stands for the Sobolev space of zero mean value functions $v$ on $S^1$ with finite norm
\begin{equation} \label{Sobolevspectr}
\left\|v\right\|_{s}=(2 \pi)^{s} \Big( \sum_{k \in \Z}{|k|^{2s} |\hat{v}(k)|^2} \Big)^{1/2},
\end{equation}
where $\hat{v}(k)$ are the complex Fourier coefficients of $v(x)$. For an integer $s=m$, this norm coincides with the previously defined $H^m$ norm. For $s \in (0,1)$, $\left\|v\right\|_{s}$ is equivalent to the norm
\begin{equation} \label{Sobolevfrac}
\left\|v\right\|^{'}_{s}=\Bigg( \int_{S^1} \Big(\int_0^1 {\frac{|v(x+\ell)-v(x)|^2}{\ell^{2s+1}} d \ell} \Big) dx \Bigg)^{1/2}
\end{equation}
(see \cite{Ada75,Tay96}).
\\ \indent
Subindices $t$ and $x$, which can be repeated, denote partial differentiation with respect to the corresponding variables. We denote by $v^{(m)}$ the $m$-th derivative of $v$ in the variable $x$. For shortness, the function $v(t,\cdot)$ is denoted by $v(t)$.

\subsection{Different types of forcing} \label{force}

In Section~\ref{rand}, we consider the generalised Burgers equation with two different types of additive forcing in the right-hand side. Since the forcing is always a r.v. in $L_2$ and the initial condition satisfies (\ref{zero}), its solutions satisfy (\ref{zero}) for all time.
\\ \indent
First, we consider the kick force. We begin by providing each space $W^{m,p}$ with the Borel $\sigma$-algebra. Then we consider an $L_2$-valued r.v. $\zeta=\zeta^{\omega}$ on a probability space $(\Omega, \mathcal{F}, \Pe)$. We suppose that $\zeta$ satisfies the following three properties.
\\ \indent
\textbf{(i) (Non-triviality)}
$$
\Pe(\zeta \equiv 0)<1.
$$
\indent
\textbf{(ii) (Finiteness of moments for Sobolev norms)}
For every
\\
$m \geq 0$, we have
$$
I_m=\E \left\|\zeta\right\|^{2}_m < +\infty,\quad \forall m \geq 0.
$$
\\ \indent
\textbf{(iii) (Vanishing of the expected value)}
$$
\E \zeta \equiv 0.
$$
\\ \indent
It is not difficult to construct explicitly
$\zeta$ satisfying \textbf{(i)-(iii)}. For instance
we could consider the real Fourier coefficients of $\zeta$, defined for $k>0$ by
\begin{equation} \label{kickFourier}
a_k(\zeta)=\sqrt{2} \int_{S^1}{\cos(2 \pi kx) u(x)};\ b_k(\zeta)=\sqrt{2} \int_{S^1}{\sin(2 \pi kx) u(x)},
\end{equation}
as independent r.v. with zero mean value and exponential moments tending to $1$ fast enough as $k \rightarrow + \infty$.
\\ \indent
Now let $\zeta_i$, $i \in \N$ be i.i.d. r.v.'s having the same distribution as $\zeta$. The sequence $(\zeta_i)_{i \geq 1}$ is a r.v. defined on a probability space which is a countable direct product of copies of  $\Omega$. From now on, this space will itself be called $\Omega$. The meaning of $\mathcal{F}$ and $\Pe$ changes accordingly.
\\ \indent
For $\omega \in \Omega$, the kick force $\xi^{\omega}$ is a $C^{\infty}$-smooth function in the variable $x$, with values in the space of distributions in the variable $t$, defined by
$$
\xi^{\omega}(x)=\sum_{i=1}^{+\infty}{\delta_{t=i} \zeta_i^{\omega}(x)},
$$
where $\delta_{t=i}$ denotes the Dirac measure at a time moment $i$.
\\ \indent
The kick-forced equation corresponds to the case where, in the right-hand side of (\ref{Burbegin}), $0$ is replaced by the kick force:
\begin{equation} \label{kickBurgers}
\frac{\partial u}{\partial t} + f'(u)\frac{\partial u}{\partial x} - \nu \frac{\partial^2 u}{\partial x^2} = \xi^{\omega}.
\end{equation}This means that  for integers $i \geq 1$, at the moments $i$ the solution $u(x)$ instantly increases by the kick $\zeta_i^{\omega}(x)$, and that between these moments $u$ solves (\ref{Burbegin}).
\\ \indent
The other type of forcing considered here is the \textit{white force}. Heuristically this force corresponds to a scaled limit of kick forces with more and more frequent kicks.
\\ \indent
To construct the white force, we begin by considering an $L_2$-valued random process 
$$
w(t)=w^{\omega}(t),\ \omega \in \Omega,\ t \geq 0,
$$
defined on a complete probability space $(\Omega, \F, \Pe)$. We assume that $w(t)$ is a Wiener process with respect to a filtration $\F_t,\ t \geq 0$, in any space $H^m,\ m \geq 0$. In particular, for $\zeta,\chi \in L_2,$
$$
\E(\left\langle w(s),\zeta\right\rangle \left\langle w(t),\chi \right\rangle)=\min(s,t) \left\langle Q\zeta,\chi\right\rangle,
$$
where $Q$ is a symmetric operator which defines a continuous mapping $Q: L_2 \rightarrow H^m$ for every $m$. Thus, $w(t) \in C^{\infty}$ for every $t$, a.s. We will denote $w(t)(x)$ by $w(t,x)$.
For more details, see \cite[Chapter 4]{DZ92}. For $m \geq 0$, we denote by $I_m$ the quantity
$$
I_m=Tr_{H^m}(Q)=\E \left\|w(1)\right\|_m^2.
$$
It is not difficult to construct $w(t)$ explicitly. For instance, we could consider the particular case of a \enquote{diagonal} noise:
\begin{equation} \nonumber
w(t)=\sqrt{2} \sum_{k \geq 1}{a_k w_k(t) \cos(2 \pi kx)}+\sqrt{2}\sum_{k \geq 1}{b_k \tilde{w}_k(t) \sin(2 \pi kx)},
\end{equation}
where $w_k(t),\ \tilde{w}_k(t),\ k > 0,$ are standard independent Wiener processes and
$$
I_m=\sum_{k \geq 1}{(a_k^2+b_k^2) (2 \pi k)^{2m}} < \infty
$$
for each $m$.  From now on, the term $dw(s)$ denotes the stochastic differential corresponding to the Wiener process $w(s)$ in the space $L_2$.
\\ \indent
Now fix $m \geq 0$. By Fernique's Theorem \cite[Theorem 3.3.1]{Kuo75}, there exist $\lambda_{m},C_m>0$ such that
\begin{equation} \label{eexp}
\E \exp \Big(\lambda_m \left\|w(T)\right\|_m^2/T \Big) \leq C_m,\quad T \geq 0.
\end{equation}
Therefore by Doob's maximal inequality for infinite-dimensional submartingales \cite[Theorem 3.8. (ii)]{DZ92} we have
\begin{equation} \label{moments}
\E \sup_{t \in [0,T]} {\left\|w(t)\right\|^p_m} \leq \Big( \frac{p}{p-1} \Big)^p \E \left\|w(T)\right\|_m^p < +\infty,
\end{equation}
for any $T>0$ and $p \in (1,\infty)$.
\\ \indent
The white-forced equation is obtained by replacing $0$ by $\eta^{\omega}=\partial w^{\omega}/\partial t$ in the right-hand side of (\ref{Burbegin}). Here, $w^{\omega}(t),\ t \geq 0$ is the Wiener process with respect to the filtration $\left\{ \mathcal{F}_t \right\}$ defined above.

\begin{defi} \label{whiteweak}
We say that an $H^1$-valued process $u(t,x)=u^{\omega}(t,x)$ is a solution of the equation
\begin{equation} \label{whiteBurgers}
\frac{\partial u^{\omega}}{\partial t} + f'(u^{\omega})\frac{\partial u^{\omega}}{\partial x} - \nu \frac{\partial^2 u^{\omega}}{\partial x^2} = \eta^{\omega}
\end{equation}
if 
\\ \indent
(i) For every $t$, $\omega \mapsto u^{\omega}(t)$ is $\mathcal{F}_t$-measurable.
\\ \indent
(ii) For a.e. (almost every) $\omega$, $t \mapsto u^{\omega}(t)$ is continuous in $H^1$ and satisfies
\begin{align} \label{whiteBurgersintbis}
u^{\omega}(t)=u^{\omega}(0)-\int_{0}^{t}{\Big( \nu Lu^{\omega}(s)+\frac{1}{2} B(u^{\omega})(s) \Big) ds}+w^{\omega}(t),
\end{align}
where
$$
B(u)=2 f'(u)u_x;\quad L=-\partial_{xx}.
$$
\end{defi}
\indent
Now consider, for a solution $u(t,x)$ of (\ref{whiteBurgers}), the functional $G_m(u(t))=\left\|u(t)\right\|_m^2$ and apply It{\^o}'s formula \cite[Theorem 4.17]{DZ92}:
\begin{align} \nonumber
\left\|u(t)\right\|_m^2=& \left\|u_0\right\|_m^2 - \int_{0}^{t}{\left( 2\nu \left\|u(s)\right\|_{m+1}^2 +\langle L^m u(s),\ B(u)(s) \rangle\right) ds}+ tI_m
\\ \nonumber
 &+ 2 \int_{0}^{t}{\langle L^m u(s),\ dw(s)\rangle}
\end{align}
(we recall that $I_m=Tr(Q_m)$.) Consequently,
\begin{align} \nonumber
\frac{d}{dt} \E \left\|u(t)\right\|_m^2 &=-2 \nu \E \left\|u(t)\right\|_{m+1}^2 - \E\ \langle L^m u(t),\ B(u)(t) \rangle+ I_m.
\end{align}
As $\langle u,\ B(u) \rangle=0$, for $m=0$ this relation becomes
\begin{align} \label{Itoexp0diff}
\frac{d}{dt} \E \left|u(t)\right|^2 &= I_0-2 \nu \E \left\|u(t)\right\|_{1}^2.
\end{align}

\subsection{Notation and agreements} \label{agree}

When considering a Sobolev norm in $W^{m,p}$, the quantity $\gamma=\gamma(m,p)$ denotes $\max(0,m-1/p)$.
\\ \indent
In Section~\ref{K41}, $\ve(t,\ix)$ denotes the velocity of a 3d flow with period $1$ in each spatial coordinate. In the whole paper, $u(t,\ix)$ denotes a solution of the generalised Burgers equation with a given initial condition $u(0,\cdot)$. In Section~\ref{det}, we deal with the equation (\ref{Burbegin}) under the assumptions (\ref{strconvex}-\ref{zero}). In Section~\ref{rand} we deal with the equation (\ref{whiteBurgers}), under the assumptions (\ref{strconvex}-\ref{zero}) and under the additional assumption
\begin{equation} \label{poly}
\forall m \geq 0,\ \exists h \geq 0,\ C_m>0:\ |f^{(m)}(x)| \leq C_m (1+|x|)^h,\quad x \in \R,
\end{equation}
where $h=h(m)$ is a function such that $1 \leq h(1) < 2$ (the lower bound on $h(1)$ follows from (\ref{strconvex})). The results in that section also hold for the kicked equation (\ref{kickBurgers}).
\\ \indent
When we consider the randomly forced generalised Burgers equation, $\Pe$ et $\E$ denote, respectively, the probability and the expected value with respect to the probability measure $\Omega$ (cf. Section~\ref{force}).
\\ \indent
All quantities denoted by $C$ with sub- or superindices are nonnegative and nonrandom. Unless otherwise stated, they only depend on the following parameters:
\begin{itemize}
\item When dealing with the K41 theory, the statistical properties of the forcing.
\\
\item When studying the unforced generalised Burgers equation, the function $f$ determining the nonlinearity $f'(u)u_x$, as well as the parameter
\begin{equation} \label{D}
D=\max ( |u_0|_{1}^{-1},\ |u_0|_{1,\infty} )
\end{equation}
which characterises how generic the initial condition is.
\\
\item When studying the randomly forced generalised Burgers equation, the function $f$ determining the nonlinearity $f'(u)u_x$, as well as the statistical properties of the forcing $\eta$. In the case of a kick force, by statistical properties we mean the distribution function of the i.i.d. r.v.'s $\eta_k$. In the case of a white force, we mean the correlation operator $Q$ for the Wiener process $w$ defining the random forcing.
\end{itemize}
In particular, those quantities never depend on the viscosity coefficient $\nu$.
\\ \indent
Constants which also depend on parameters $a_1,\dots,a_k$ are denoted by $C(a_1,\dots,a_k)$. By $X \overset{a_1,\dots,a_k}{\lesssim} Y$ we mean that $X \leq C(a_1,\dots,a_k) Y$. The notation $X \overset{a_1,\dots,a_k}{\sim} Y$ stands for
$$
Y \overset{a_1,\dots,a_k}{\lesssim} X \overset{a_1,\dots,a_k}{\lesssim} Y.
$$
In particular, $X \lesssim Y$ and $X \sim Y$ mean that $X \leq C Y$ and $C^{-1} Y \leq X \leq C Y$, respectively.
\\ \indent
The initial condition $u(0,\cdot)$ is denoted by $u_0$.
\\ \indent
We use the notation $g^{-}=\max(-g,0)$ and $g^{+}=\max(g,0)$.
\\ \indent
In Subsection~\ref{K41}, the brackets $\langle \cdot \rangle$ denote the expected value. For the meaning of the brackets $\lbrace \cdot \rbrace$, see Subsection~\ref{detSob} in the deterministic case (where they correspond to averaging in time) and Subsection~\ref{whiteturb} in the random case (where they correspond to averaging in time and taking the expected value). The definitions of the small-scale quantities, i.e. the structure functions $S_{p,\alpha}$ and $S_{p,1}=S_p$ and the spectrum $E(k)$ depend on the setting: see Subsections~\ref{K41}, \ref{Burgu}, \ref{detturb} and \ref{whiteturb}.

\section{Turbulence and the Burgers equation} \label{turb}

\subsection{Turbulence, K41 theory, intermittency} \label{K41}

\begin{figure}[t]
\includegraphics[width=4.4cm,height=2.8cm]{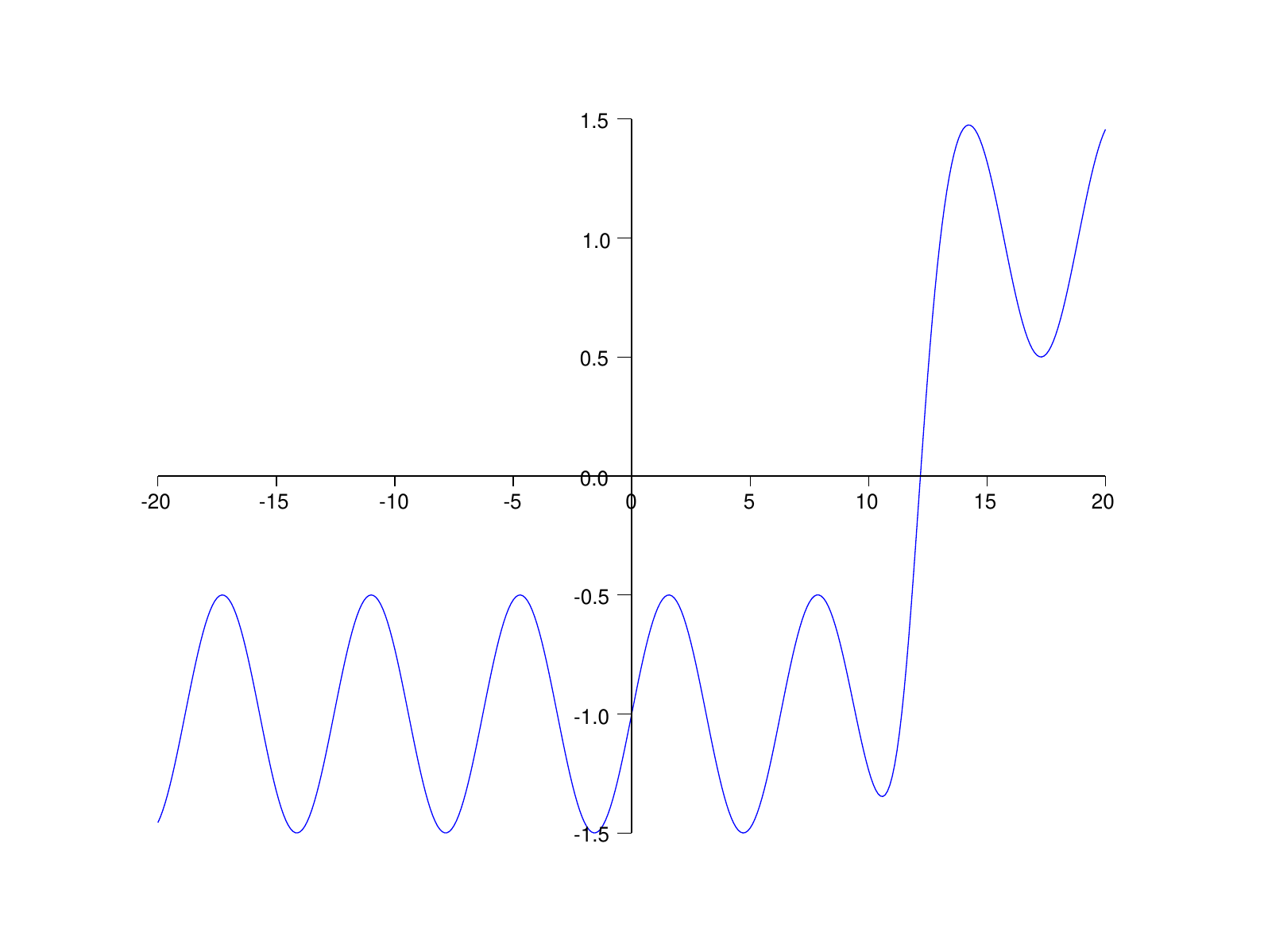}
\caption{Example of a function exhibiting small-scale intermittency}
\label{intergraph}
\end{figure}

It is well-known that giving a precise definition of turbulence is problematic. However, some features are generally recognised as characteristic of turbulence: presence of many degrees of freedom, unpredictability/chaos, (small-scale) irregularity... For a more detailed discussion, see \cite{Fri95,Tsi09}. Here, we will only present (in a slightly modified form) the vocabulary of the theory of turbulence which is relevant to the study of the Burgers model. In particular, we will proceed as if the flow $\ve(t,\ix)$ is periodic in space. Without loss of generality, we may assume that $\ve$ is $1$-periodic in each coordinate $x_1,x_2,x_3$. Let us denote by $\nu$ the viscosity coefficient; we only consider the turbulent regime $0 < \nu \ll 1$.
\\ \indent
We define the \textit{space scale} as the inverse of the frequency under consideration. In particular, the Fourier coefficients $\hat{\ve}(\ka)$ for large values of $\ka$ or, in the physical space, the increments $\ve(\ix+\er)-\ve(\ix)$ for small values of $\er$, are prototypical small-scale quantities.
\\ \indent
The theory which may be considered as a starting point for the modern study of turbulence is essentially contained in three articles by Kolmogorov which have been published in 1941 \cite{Kol41a,Kol41b,Kol41c}. Thus, it is referred to as the \textit{K41 theory}.
\\ \indent
The philosophy behind K41 is that although large-scale characteristics of a turbulent flow are clearly \enquote{individual} (depending on the forcing or on the boundary conditions), small-scale characteristics display some non-trivial \enquote{universal} features. To make this point clearer, we will introduce several definitions.
\\ \indent
The \textit{dissipation scale} $\ell_d$ is the smallest scale such that for all $| \ka | \succeq \ell_d^{-1}$, the Fourier coefficients of a function $\ve$ decrease super-polynomially in $| \ka |$, uniformly in $\nu$. The interval $\J_{diss}=(0,\ell_d]$ is called the  \textit{dissipation range}. The K41 theory claims that $\ell_d=C \nu^{3/4}$. The \textit{energy range} $\J_{energ}=(\ell_e,1]$ consists of the scales such that the corresponding Fourier modes support most of the $L^2$ norm of $\ve$:
$$
\sum_{| \ka | < \ell_e^{-1}}{\langle|\hat{\ve}(\ka)|^2\rangle} \gg \sum_{| \ka | \geq \ell_e^{-1}}{\langle|\hat{\ve}(\ka)|^2\rangle}.
$$
K41 states that $\ell_e=C$.
\\ \indent
Finally $\J_{inert}=(\ell_d,\ell_e]$ is the \textit{inertial range}. K41 states that $\J_{inert}=(C \nu^{3/4},C]$. This is the most interesting zone, where the flow exhibits non-trivial small-scale behaviour which will be described more precisely below.
\smallskip
\begin{figure}[h]
\includegraphics[height=2cm]{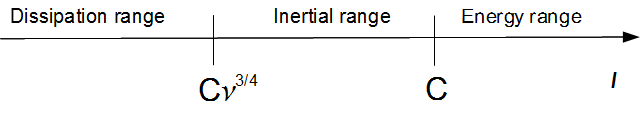}
\caption{Kolmogorov scales}
\end{figure}
\smallskip \indent \\
Two quantities used to describe small-scale behaviour of a flow $\ve(t,\ix)$ \textit{at a fixed time moment} are:
\begin{itemize}
\item On one hand, the \textit{longitudinal structure function}
\begin{equation} \label{introSp}
S^{\Vert}_p(\ix,\er)= \Bigg< \Bigg| \frac{(\ve(\ix+\er)-\ve(\ix)) \cdot \er}{|\er|} \Bigg|^p \Bigg>
\end{equation}
\item On the other hand, the \textit{energy spectrum} 
\begin{equation} \label{introEk}
E(k)=\frac{\sum_{|\en| \in [M^{-1}k,Mk]}{\langle |\hat{\ve}(\en)|^2 \rangle}}{\sum_{|\en| \in [M^{-1}k,Mk]}{1}},
\end{equation}
i.e. the average of $\langle |\hat{\ve}(\en)|^2 \rangle$ over a layer of $\en$ such that $|\en| \sim k$.
\end{itemize}
\indent
The K41 theory predicts that under some conditions on the flow, for $\ell=|\er| \in \J_{inert}$ and for every $\ix$, we have
\begin{equation} \label{introK41phy}
S^{\Vert}_p(\ix,\er) \overset{p}{\sim} \ell^{p/3},\quad p \geq 0.
\end{equation}
On the other hand, for $k$ such that $k^{-1} \in \J_{inert}$, K41 states that
\begin{equation} \label{introK41Fou}
E(k) \sim k^{-5/3}
\end{equation}
(see \cite{Obu41a,Obu41b}).
\\ \indent
The K41 predictions are in good agreement with experimental and numerical data for the energy spectrum and for the structure functions $S_p,\ p=2,3$. However, there are important discrepancies for the functions $S_p,\ p \geq 4$ \cite[Chapter 8]{Fri95}. Two parallel theories, due respectively to Kolmogorov himself \cite{Kol62} and to Frisch and Parisi \cite{ParFri85} give an explanation which emphasises the role of \textit{spatial intermittency}. In other words, at a given time moment, the flow is very strongly excited on a small subset, as for the function whose graph is given in Figure~\ref{intergraph}.
\\ \indent
Intermittency at the scale $\ell$ is quantified by \textit{flatness}, defined as
$$
F(\ell)=S^{\Vert}_4(\ell)/S^{\Vert}_2(\ell)^2:
$$
the larger the flatness, the more intermittent is the function. Thus, we need to take into account the intermittency since the K41 theory does not predict the corresponding features observed in the inertial range in turbulent flows such as vortex stretching \cite{SheOrs91}: indeed, for $\ell \in J_2$ the K41 predictions yield that
$$
F(\ell) \sim \ell^{4/3}/(\ell^{2/3})^2 =1.
$$

\subsection{Burgers turbulence} \label{Burgu}
\indent
The 1d Burgers equation
\begin{equation} \label{introBur}
u_t+uu_x=\nu u_{xx},
\end{equation}
where $\nu>0$ is a viscosity coefficient, has first been considered by Forsyth \cite{For06} and Bateman \cite{Bate15} in the first decades of the XXth century. Here, we will only consider the space-periodic case: after rescaling, we can suppose that $\ix \in S^1=\R/\Z$.
\\ \indent
This equation is well-posed in $L_1(S^1)$. Indeed, the proof of such a statement in a smaller space is very standard: see for instance \cite[Chapter 5]{KrLo89}. Well-posedness in $L_1(S^1)$ follows then by a contraction argument (see Section~\ref{stat}).
\\ \indent
Around 1950, the Burgers equation attracted considerable interest in the scientific community. In particular, it has been studied by the Dutch physicist whose name it bears (\cite{Bur48, Bur74}; see also \cite{Batc53}). His goal was to consider a simplified version of the incompressible Navier--Stokes equation
\begin{equation} \label{introNS}
\ue_t+(\ue \cdot \nabla) \ue=\nu \Delta \ue-\nabla p;\quad \nabla \cdot \ue=0,
\end{equation}
which would keep some of its features. This hope was shared by von Neumann \cite[p. 437]{Neu63}. 
\\ \indent
The Hopf-Cole-Florin transformation (\cite{Col51,Flo48,Hop50}; see \cite{Bir03} for a historical account) reduces the Burgers equation to the heat equation. Indeed, if $u$ is the solution of (\ref{introBur}) corresponding to an initial condition $u_0$, then $u(t,x)$ is the space derivative of the function
$$
-2 \nu \ln(\phi(t,x)),
$$
where $\phi$ is the solution of the heat equation
$$
\phi_t=\nu \phi_{xx}
$$
corresponding to the initial condition $\phi_0=\exp(-H_0/2 \nu)$. Here, $H_0$ is a primitive of $u_0$. This transformation can also be applied to the multi-d potential Burgers equation:
\begin{equation} \label{intromultiBur}
\ue_t+(\ue \cdot \nabla) \ue=\nu \Delta \ue;\quad \ue=-\nabla \psi.
\end{equation}
\indent
Note that such a transformation does not exist for the generalised Burgers equation considered in our survey.
\\ \indent
The fact that the Burgers equation can be reduced to the heat equation means that it is integrable and therefore its solutions do not exhibit chaotic behaviour. However, the Hopf-Cole-Florin transformation cannot immediately provide information about the small-scale behaviour of solutions in the turbulent regime corresponding to $0 < \nu \ll 1$. This behaviour has been studied on a qualitative level by many physicists \cite{AFLV92,Cho75,Kra68,Kid79}. There is an agreement about the behaviour of the increments and of the energy spectrum in the inertial range, which corresponds to the interval $\J_{inert}=(C \nu, C]$.
\smallskip
\begin{figure}[h]
\includegraphics[height=2cm]{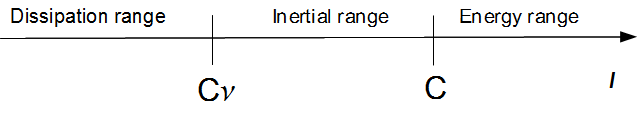}
\caption{Space scales for the Burgers equation}
\end{figure}
\smallskip \indent \\
First, if we denote by $S_p(\ell)$ the structure function defined by
\begin{equation} \label{defSp}
S_p(\ell)= \int_{S^1}{ | u(x+\ell)-u(x) |^p dx },
\end{equation}
then for $\ell \in \J_{inert}$ we have
\begin{equation} \label{valSp}
S_p(\ell) \overset{p}{\sim} \left\lbrace \begin{aligned} & \ell^p,\ 0 \leq p \leq 1. \\ &\ell,\ p \geq 1. \end{aligned} \right.
\end{equation}
In particular, for $\ell$ in the inertial range, the flatness $F(\ell)$ behaves as $\ell^{-1}$. This is related to the intermittent behaviour on small scales corresponding to the \enquote{cliffs} of a typical solution, which will be described below.
\\ \indent
On the other hand, for $k^{-1} \in \J_{inert}$ we have $E(k) \sim k^{-2}$ with the same definition as above (up to the absence of the brackets $\left\langle \cdot \right\rangle$) for $E(k)$.
\\ \indent
To explain the physical arguments of \cite{AFLV92}, we need to give more details on the structure of solutions for (\ref{introBur}). We assume that both the initial condition $u_0$ and its derivative have amplitude of the order 1.
\\ \indent
First, consider the inviscid Hopf equation which is the limit case $\nu=0$ of (\ref{introBur}). Its solution is only smooth during a finite interval of time: it can be implicitly constructed using the method of characteristics (see for instance \cite{Daf10}). This method tells us that while the solution remains smooth, the value of $u$ is constant along the lines $(t,x+tu_0(x))$ in the space-time. However, if $u_0$ is not constant, then lines corresponding to different values of $u_0$ cross after a finite time, forbidding the existence of smooth solutions. Nevertheless, a weak entropy solution can still be uniquely defined for all time in the class of bounded variation functions $BV(S^1)$. Such a solution is a limit in $L_1$ of classical solutions for the viscous equation as $\nu \rightarrow 0$. More precisely, this solution exhibits the $N$-wave behaviour \cite{Eva08}, i.e. for a fixed time $t$ its graph is similar to repeated mirror images of the capital letter N. In other words, the solution $u(t,\cdot)$ alternates between negative jump discontinuities and smooth regions where the derivative is positive and of order $1$.
\\ \indent
When $\nu>0$, the shocks become cliffs. The amplitude of the solution, the number of cliffs and the height of a cliff are all of order $1$. The width of a cliff is of order $\nu$.
\begin{figure}[h]
\includegraphics[width=9 cm, height=5cm]{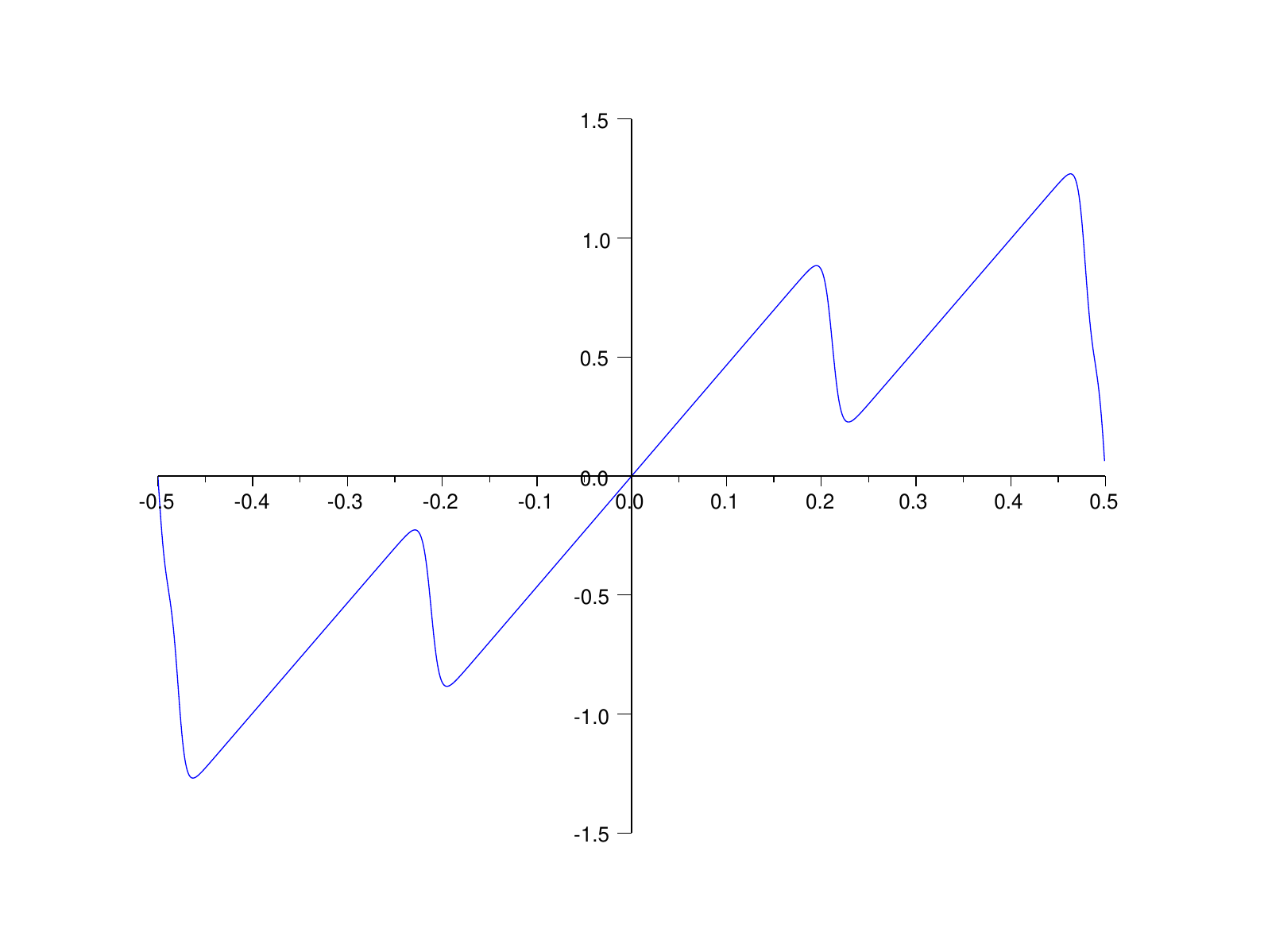}
\caption{\enquote{Typical} solution of the Burgers equation}
\end{figure}
\\ \indent
For $\ell \in \J_{inert}$, $\ell$ is typically smaller than the interval between two cliffs, but larger than the width of a cliff. Aurell, Frisch, Lutsko and Vergassola observe that there are $3$ possibilities for the interval $[x,x+\ell]$.
\smallskip
\begin{itemize}
\item{$[x,x+\ell]$ covers a large part of a "cliff".}
\\
Probability $=$ \textit{$C \ell$.}\quad $u(x+\ell)-u(x)=\underbrace{-C}_{cliff}+\underbrace{C \ell}_{ramps}=-C$.\quad $|u(x+\ell)-u(x) |^p \overset{p}{\sim} \textit{1.}$
\\
\item{$[x,x+\ell]$ covers a small part of a "cliff".}
\\
Contribution of this term is negligible.
\\
\item{$[x,x+\ell]$ does not intersect a "cliff".}
\\
Probability $= 1-C \ell = C$.\quad $u(x+\ell)-u(x)=\underbrace{C \ell}_{ramp}$.
\\
$|u(x+\ell)-u(x) |^p \overset{p}{\sim} \ell^p.$
\end{itemize}
\smallskip
\indent
Thus, $S_p(\ell) \overset{p}{\sim} \ell+\ell^p \overset{p}{\sim} \left\lbrace \begin{aligned} & \ell^p,\ 0 \leq p \leq 1. \\ &\ell,\ p \geq 1. \end{aligned} \right.$
\\ \indent
In other words, for $p \geq 0$ the description above implies that for $\ell \in \J_{inert}$, the behavour of the structure functions is given by (\ref{valSp}). 
\\ \indent
Asymptotically, the Fourier coefficients of an $N$-wave satisfy $|\hat{u}(k)| \sim k^{-1}$. Thus, it is natural to conjecture that for $\nu$ small and for a certain range of $k$, energy-type quantities $\frac{1}{2} |\hat{u}(k)|^2$ behave, in average, as $k^{-2}$ \cite{Cho75,FouFri83,Kid79,Kra68}.
\\ \indent
Beginning from the 1980s, there has been an increasing interest in random versions of the Burgers equation. The most studied model has been the one with additive white in time noise, more or less smooth in space. Here, we will only consider the case where the noise is
\\
$C^{\infty}$-smooth in space; for the general case, see the surveys \cite{BF01,BK07}. In that setting, numerical simulations and physical predictions give exactly the same results as in the deterministic case, up to the fact that we consider the expected values of the quantities \cite{GotKra98}. Heuristically, this is due to the fact that forcing acts on large scales, in the energy range, and thus only influences smaller scales indirectly, as an energy source.

\section{Main results} \label{results}

In Section~\ref{det}, we are concerned with the deterministic Burgers
\\
equation. First, in Subsection~\ref{detSob}, we prove sharp upper and lower bounds for some Sobolev norms of $u$. In Lemma~\ref{detuxpos}, we recall the key estimate
\begin{equation}\label{detuxposresults}
u_x(t,x) \leq \min (D,\sigma^{-1} t^{-1}).
\end{equation}
The main results for Sobolev norms of solutions are summed up in Theorem~\ref{detavoir}. Namely, for $m=0,1$ and $p \in [1,\infty]$ or for $m \geq 2$ and  $p \in (1,\infty]$, we have
\begin{equation} \label{detavoirresults}
\Big( \lbrace {\left|u(t)\right|_{m,p}^{\alpha}} \rbrace \Big)^{1/\alpha} \overset{m,p,\alpha}{\sim} \nu^{-\gamma},\quad \alpha>0,
\end{equation}
where $\lbrace \cdot \rbrace$ denotes averaging in time over the interval $[T_1,T_2]$ defined by (\ref{detT1T2}). We recall that $\gamma(m,p)=\max(0,m-1/p)$.
\\ \indent
In Subsection~\ref{detturb} we obtain sharp estimates for analogues of the quantities characterising the hydrodynamical turbulence. In what follows, we assume that $\nu \in (0,\nu_0]$, where $\nu_0 \in (0,1]$ only depends on $f$ and on $D$. To begin with, we define the non-empty and non-intersecting intervals
\begin{equation} \nonumber
J_1=(0,\ C_1 \nu];\ J_2=(C_1 \nu,\ C_2];\ J_3=(C_2,\ 1].
\end{equation}
For the definitions of $\nu_0$, $C_1$ and $C_2$, see (\ref{detC1C2}); those quantities only depend on $f$ and on $D$. As a consequence of (\ref{detuxposresults}-\ref{detavoirresults}), in Theorem~\ref{detavoir2} we prove that for $\ell \in J_1$:
$$
\quad \ \ \ S_{p}(\ell) \overset{p}{\sim} \left\lbrace \begin{aligned} & \ell^{p},\ 0 \leq p \leq 1. \\ & \ell^{p} \nu^{-(p-1)},\ p \geq 1, \end{aligned} \right.$$
and for $\ell \in J_2$:
$$
S_{p}(\ell) \overset{p}{\sim} \left\lbrace \begin{aligned} & \ell^{p},\ 0 \leq p \leq 1. \\ & \ell,\ p \geq 1. \end{aligned} \right.
$$
Consequently, for $\ell \in J_2$ the flatness satisfies:
$$
F(\ell):=S_4(\ell)/S_2^2(\ell) \sim \ell^{-1}.
$$
Finally, we get estimates for the spectral asymptotics of Burgulence. On one hand, as a consequence of Theorem~\ref{detavoir}, for $m \geq 1$ we get:
$$
\lbrace |\hat{u}(k)|^2 \rbrace \overset{m}{\lesssim} k^{-2m} {\Vert u \Vert_m^2 } \overset{m}{\lesssim} (k \nu)^{-2m} \nu.
$$
In particular, $\lbrace |\hat{u}(k)|^2 \rbrace$ decreases at a faster-than-algebraic rate for $|k| \succeq \nu^{-1}$. On the other hand, by Theorem~\ref{detspectrinert}, for $k$ such that $k^{-1} \in J_2$ the energy spectrum satisfies
$$
E(k)=\frac{\sum_{|n| \in [M^{-1}k,Mk]}{\langle |\hat{u}(n)|^2 \rangle}}{\sum_{|n| \in [M^{-1}k,Mk]}{1}} \sim  k^{-2},
$$
where $M \geq 1$ depends only on $f$ and on $D$.
\\ \indent
Note that these results rigorously confirm the physical predictions exposed in Subsection~\ref{Burgu}. Moreover, averaging in the initial condition, as considered in \cite{AFLV92}, is actually not necessary. This is due to the particular structure of the deterministic generalised Burgers equation: an initial condition $u_0$ is as \enquote{generic} as the ratio between the orders of $(u_0)_x$ and of $u_0$ itself, which can be bounded from above using the quantity $D$.
\\ \indent
The results in Section~\ref{rand} can be formulated in exactly the same way, up to three modifications:
\begin{itemize}
\item All quantities should be replaced by their expected values. In particular, we modify the meaning of the brackets $\lbrace \cdot \rbrace $.
\item Dependence on $D$ should be replaced by dependence on the statistical properties of the forcing.
\item The estimates hold uniformy in $t$ (for $t$ large enough) and in $u_0$.
\end{itemize}
\indent
In Section~\ref{stat}, we expose results on existence and uniqueness of the stationary measure for the randomly forced generalised Burgers equation. These results yield that all estimates listed above still hold with taking the expected value and averaging in time replaced by averaging with respect to the stationary measure $\mu$.

\section{The deterministic Burgers equation} \label{det}

The results in Subsection~\ref{detSob} have been obtained in \cite{Bir01} for norms in $H^m,\ m \geq 1$, under a slightly different form. Our presentation follows the lines of \cite{BorD}, where some additional estimates on Sobolev norms are obtained by H{\"o}lder's inequality and (GN). In \cite{Bir01}, Biryuk also proved upper and lower spectral estimates. The sharp small-scale results in Subsection~\ref{detturb} have been obtained in \cite{BorD}.

\subsection{Estimates for Sobolev norms} \label{detSob}

We begin by recalling the proof of a key upper estimate for $u_x$, which is a reformulation of the 
\\
\enquote{Kruzhkov maximum principle} \cite{Kru64}.

\begin{lemm} \label{detuxpos}
We have
$$
u_x(t,x) \leq \min (D,\sigma^{-1} t^{-1} ).
$$
\end{lemm}

\textbf{Proof.}
Differentiating the equation (\ref{Burbegin}) once in space we get
\begin{equation} \label{uxdiff}
(u_x)_t+f''(u)u_x^2+f'(u)(u_{x})_{x}=\nu (u_x)_{xx}.
\end{equation}
Now consider a point $(t_1,x_1)$ where $u_x$ reaches its maximum on the cylinder $S=[0,t] \times S^1$. Suppose that $t_1>0$ and that this maximum is nonnegative. At such a point, Taylor's formula implies that we would have $(u_x)_t \geq 0$, $(u_{x})_{x}=0$ and $(u_x)_{xx} \leq 0$. Consequently, since by (\ref{strconvex}) $f''(u) \geq \sigma$, (\ref{uxdiff}) yields that $\sigma u_x^2 \leq 0$, which is impossible. Thus $u_x$ can only reach a nonnegative maximum on $S$ for $t_1=0$. In other words, since $(u_0)_x$ has zero mean value, we have
$$
u_x(t,x) \leq \max_{x \in S^1} {(u_0)_x(x)} \leq D.
$$
\\ \indent
The inequality
$$
u_x(t,x) \leq \sigma^{-1} t^{-1}
$$
is proved in \cite{Kru64} by a similar maximum principle argument applied to the function $v=tu_x$. Indeed, this function can only reach a nonnegative maximum on $S$ at a point $(t_1,x_1)$ such that $t_1>0$. Multiplying (\ref{uxdiff}) by $t^2$, we get: 
$$
t \underbrace{v_t}_{\textrm{$\geq 0$}}+t f'(u) \underbrace{v_x}_{\textrm{0}}+
(-v+f''(u)v^2)=\nu t \underbrace{v_{xx}}_{\textrm{$\leq 0$}}.
$$
Thus $v \leq \sigma^{-1}$ on $S$. In other words, $u_x \leq \sigma^{-1} t^{-1}$ for all $t>0$.\ $\square$
\smallskip
\\ \indent
Since the space averages of $u(t)$ and $u_x(t)$ vanish for all $t$, we get the following upper estimates:
\begin{align} \label{detLpupper}
&\left|u(t)\right|_{p} \leq \left|u(t)\right|_{\infty} \leq \int_{S^1}{u_x^{+}(t)} \leq \min (D, \sigma^{-1} t^{-1} ),\quad 1 \leq p \leq +\infty.
\\ \label{detu11}
&\left|u(t)\right|_{1,1}=2 \int_{S^1}{u_x^{+}(t)} \leq 2\min (D, \sigma^{-1} t^{-1}).
\end{align}

Now we recall a standard estimate for the nonlinearity
$$
\left\langle v^{(m+1)}, (f(v))^{(m)}\right\rangle,
$$
which is proved in \cite{BorW}.

\begin{lemm} \label{detlmubuinfty}
For $v \in C^{\infty}$ such that $\left|v\right|_{\infty} \leq N$, we have
$$
N_m(v)=\left| \left\langle  v^{(m+1)}, (f(v))^{(m)} \right\rangle\right| \overset{m,N}{\lesssim} \left\|v\right\|_m \left\|v\right\|_{m+1},\quad m \geq 1.
$$
\end{lemm}

\textbf{Proof.}
Fix $m \geq 1$. In this proof, constants denoted by $\tilde{C}$ only depend on $m,N$. We have
\begin{align} \nonumber
&N_m(v) \leq  \tilde{C} \sum_{k=1}^m\ \sum_{\substack{1 \leq a_1 \leq \dots \leq a_k \leq m \\ a_1+ \dots+a_k = m}} \int_{S^1}{\left| v^{(m+1)} v^{(a_1)} \dots v^{(a_k)} f^{(k)}(v) \right|}
\\ \nonumber
&\leq \tilde{C} \max_{x \in [-N,N]}\ \max(f'(x),\dots f^{(m)}(x))
\\ \nonumber
& \times \sum_{k=1}^m\ \sum_{\substack{1 \leq a_1 \leq \dots \leq a_k \leq m \\ a_1+ \dots+a_k = m}} \int_{S^1} | v^{(a_1)} \dots v^{(a_k)} v^{(m+1)} |.
\end{align}
Using (\ref{poly}), H{\"o}lder's inequality and (GN), we get
\begin{align} \nonumber
&N_m(v) \leq  \tilde{C} (1+N)^{\max(h(1),\dots,h(m))} 
\\ \nonumber
&\times \sum_{k=1}^m\ \sum_{\substack{1 \leq a_1 \leq \dots \leq a_k \leq m \\ a_1+ \dots+a_k = m}} \int_{S^1} | v^{(a_1)} \dots v^{(a_k)} v^{(m+1)} |
\\ \nonumber
&\leq  \tilde{C} \sum_{k=1}^m\ \sum_{\substack{1 \leq a_1 \leq \dots \leq a_k \leq m \\ a_1+ \dots+a_k = m}} ( \left|v^{(a_1)}\right|_{2m/a_1} \dots \left|v^{(a_k)}\right|_{2m/a_k} \left\|v\right\|_{m+1} )
\\ \nonumber
&\leq  \tilde{C} \left\|v\right\|_{m+1} \sum_{k=1}^m\ \sum_{\substack{1 \leq a_1 \leq \dots \leq a_k \leq m \\ a_1+ \dots+a_k = m}} \Big( (\left\|v\right\|_m^{a_1/m} |v|_{\infty}^{(m-a_1)/m} ) \times \dots
\\ \nonumber
& \dots \times (\left\|v\right\|_m^{a_k/m} |v|_{\infty}^{(m-a_k)/m}) \Big)
\\ \nonumber
&\leq \tilde{C} (1+N)^{m-1} \left\|v\right\|_m \left\|v\right\|_{m+1} = \tilde{C} \left\|v\right\|_m \left\|v\right\|_{m+1}.\ \qed
\end{align}
\medskip \\ \indent
The following result shows the existence of a strong nonlinear damping which prevents the successive derivatives of $u$ from becoming too large.

\begin{lemm} \label{detuppermaux}
We have
$$
\left\|u(t)\right\|^{2}_1 \lesssim \nu^{-1}.
$$
On the other hand, for $m \geq 2$,
$$
\left\|u(t)\right\|^{2}_m \overset{m}{\lesssim} \max (\nu^{-(2m-1)}, t^{-(2m-1)}).
$$
\end{lemm}

\textbf{Proof.} Fix $m \geq 1$. Denote
$$
x(t)=\left\|u(t)\right\|^{2}_m.
$$
We claim that the following implication holds:
\begin{align} \label{detdecrm}
&x(t) \geq C' \nu^{-(2m-1)} \Longrightarrow \frac{d}{dt} x(t) \leq -(2m-1) x(t)^{2m/(2m-1)},
\end{align}
where $C'$ is a fixed nonnegative number, chosen later. Below, all constants denoted by $C$ do not depend on $C'$.
\\ \indent
Indeed, assume that $ x(t) \geq C' \nu^{-(2m-1)}.$ Integrating by parts in space and using (\ref{detLpupper}) ($p=\infty$) and Lemma~\ref{detlmubuinfty}, we get the following energy dissipation relation:
\begin{align} \nonumber
\frac{d}{dt}  x(t) &= - 2 \nu  \left\|u(t)\right\|_{m+1}^2+2\left\langle u^{(m+1)}(t), (f(u(t)))^{(m)}\right\rangle
\\ \indent
&\leq - 2 \nu  \left\|u(t)\right\|_{m+1}^2 + C \left\|u(t)\right\|_m \left\|u(t)\right\|_{m+1}.
\end{align}
Applying (GN) to $u_x$ and then using (\ref{detu11}), we get
\begin{align} \nonumber
\left\| u(t)\right\|_m &\leq C \left\| u(t)\right\|_{m+1}^{(2m-1)/(2m+1)} \left| u(t)\right|_{1,1}^{2/(2m+1)} 
\\ \label{detGN11m}
&\leq C \left\| u(t)\right\|_{m+1}^{(2m-1)/(2m+1)}.
\end{align}
Thus, we have the relation
\begin{align} \label{detint1}
\frac{d}{dt} x(t) \leq & (- 2 \nu \left\|u(t)\right\|_{m+1}^{2/(2m+1)} + C ) \left\|u(t)\right\|_{m+1}^{4m/(2m+1)}.
\end{align}
The inequality (\ref{detGN11m}) yields that
\begin{equation} \label{detint2}
\left\| u(t)\right\|_{m+1}^{2/(2m+1)} \geq C x(t)^{1/(2m-1)},
\end{equation}
and then since by assumption $ x(t) \geq C' \nu^{-(2m-1)}$ we get
\begin{align} \label{detint3}
 \left\|u(t)\right\|_{m+1}^{2/(2m+1)} & \geq C C'^{1/(2m-1)} \nu^{-1}.
\end{align}
Combining the inequalities (\ref{detint1}-\ref{detint3}), for $C'$ large enough we get
\begin{align} \nonumber
\frac{d}{dt}  x(t) &\leq (-C C'^{1/(2m-1)} + C ) x(t)^{2m/(2m-1)}.
\end{align}
Thus we can choose $C'$ in such a way that the implication (\ref{detdecrm}) holds.
\\ \indent
For $m=1$, (\ref{D}) and (\ref{detdecrm}) yield that
$$
x(t) \leq \max(C'\nu^{-1},D^2) \leq \max(C',D^2) \nu^{-1},\ t \geq 0.
$$
Now consider the case $m \geq 2$. We claim that
\begin{equation} \label{detdecrmcor}
x(t) \leq \max ( C' \nu^{-(2m-1)}, t^{-(2m-1)} ).
\end{equation}
Indeed, if $x(s) \leq C' \nu^{-(2m-1)}$ for some $s \in \left[0,t\right]$, then the assertion (\ref{detdecrm}) ensures that $x(s)$ remains below this threshold up to time $t$.
\\ \indent
Now, assume that $x(s) > C' \nu^{-(2m-1)}$ for all $s \in \left[0,t\right]$. Denote
$$
\tilde{x}(s)=(x(s))^{-1/(2m-1)},\ s \in \left[0,t\right].
$$
By (\ref{detdecrm}) we get $d\tilde{x}(s)/ds \geq 1$. Therefore $\tilde{x}(t) \geq t$ and $x(t) \leq t^{-(2m-1)}$. Thus in this case, the inequality (\ref{detdecrmcor}) still holds. This proves the lemma's assertion. $\square$
\smallskip
\\ \indent
By (GN) applied to $u^{(m)}$ we get the following inequality for $m \geq 1$:
\begin{align} \nonumber
\left|u(t)\right|_{m,\infty} & \lesssim \left\|u(t)\right\|^{1/2}_m \left\|u(t)\right\|^{1/2}_{m+1} \overset{m}{\lesssim} \max(\nu^{-m},t^{-m}).
\end{align}
Similarly, applying (GN) and interpolating between $|u|_{1,1}$ and $\Vert u \Vert_M$ for large values of $M$, we get the following result (we recall that $\gamma=\max (0,m-1/p)$):

\begin{theo} \label{upperwmp}
For $m \in \lbrace 0,1 \rbrace$ and $p \in [1,\infty]$, or for $m \geq 2$ and $p \in (1,\infty]$,
\begin{align} \label{detupperwmp}
\Big( \E  \max_{s \in [t,t+1]} \left|u(s)\right|^{\alpha}_{m,p} \Big)^{1/\alpha} \overset{m,p,\alpha}{\lesssim} \max(t^{-\gamma},\nu^{-\gamma}),\quad \alpha>0.
\end{align}
\end{theo}

Now we define
\begin{equation} \label{detT1T2}
T_1=\frac{1}{4}D^{-2} \tilde{C}^{-1};\quad T_2=\max \Big( \frac{3}{2} T_1,\quad 2D\sigma^{-1} \Big),
\end{equation}
where $\tilde{C}$ is a constant such that for all $t$, $\left\|u(t)\right\|_1^2 \leq \tilde{C} \nu^{-1}$ (cf.
\\
Lemma~\ref{detuppermaux}). From now on, for any function $A(t)$, $\lbrace A(t) \rbrace$ is by definition the time average
$$
\lbrace A(t) \rbrace=\frac{1}{T_2-T_1} \int_{T_1}^{T_2}{A(t)}.
$$

\begin{lemm} \label{detfinitetime}
We have
$$
\lbrace \left\|u(t)\right\|_1^2 \rbrace \gtrsim \nu^{-1}.
$$
\end{lemm}

\textbf{Proof.} Integrating by parts in space, we get the dissipation identity
\begin{align} \label{detdissip}
\frac{d}{dt} \left| u(t)\right|^2 &= \int_{S^1}{(-2 uf'(u)u_x+2 \nu u u_{xx})}=-2 \nu \left\| u(t)\right\|_{1}^{2}.
\end{align}
Thus, integrating in time and using (\ref{D}) and Lemma~\ref{detuppermaux}, we obtain that
\begin{align} \nonumber
|u(T_1)|^2 &= |u_0|^2 - 2 \nu \int_{0}^{T_1}{\left\|u(t)\right\|_1^2} \geq D^{-2} - 2 T_1 \tilde{C} \geq \frac{1}{2}D^{-2}.
\end{align}
Consequently, integrating (\ref{detdissip}) in time and using (\ref{detLpupper}) ($p=2$) we get
\begin{align} \nonumber
\lbrace \left\|u(t)\right\|_1^2 \rbrace &= \frac{1}{2 \nu (T_2-T_1)} (|u(T_1)|^2-|u(T_2)|^2)
\\ \nonumber
& \geq \frac{1}{2 \nu (T_2-T_1)} \Big(\frac{1}{2}D^{-2}-\sigma^{-2} T_2^{-2} \Big)
\\ \nonumber
& \geq \frac{D^{-2}}{8 (T_2-T_1) } \nu^{-1},
\end{align}
which proves the lemma's assertion.\ $\square$
\medskip \\ \indent
This time-averaged lower bound yields similar bounds for other
\\
Sobolev norms.

\begin{lemm} \label{detfinalexp}
For $m \geq 1$,
$$
\lbrace \left\|u(t)\right\|_m^2 \rbrace \overset{m}{\gtrsim} \nu^{-(2m-1)}.
$$
\end{lemm} 

\textbf{Proof.}
Since the case $m=1$ has been treated in the previous lemma, we may assume that $m \geq 2$. By (\ref{detu11}) and (GN), we have:
\begin{align} \nonumber
\lbrace \left\|u(t)\right\|_m^2 &\rbrace \overset{m}{\gtrsim} \lbrace \left\|u(t)\right\|_m^2  \left|u(t)\right|_{1,1}^{(4m-4)} \rbrace \overset{m}{\gtrsim} \lbrace \left\|u(t)\right\|^{4m-2}_1 \rbrace.
\end{align}
Thus, using H{\"o}lder's inequality and Lemma~\ref{detfinitetime}, we get:
\begin{align} \nonumber
\lbrace \left\|u(t)\right\|_m^2 &\rbrace \overset{m}{\gtrsim} \lbrace \left\|u(t)\right\|^{4m-2}_1 \rbrace \overset{m}{\gtrsim} \lbrace \left\|u(t)\right\|^{2}_1 \rbrace^{(2m-1)} \overset{m}{\gtrsim} \nu^{-(2m-1)}.\ \square
\end{align}
\\ \indent
\smallskip
The following lemma is proved similarly.

\begin{lemm} \label{detfinalexpbis}
For $m \geq 0$, $p \in [1,\infty]$, we have:
$$
\lbrace \left|u(t)\right|_{m,p}^2 \rbrace \overset{m,p}{\gtrsim} \nu^{-2\gamma}.
$$
\end{lemm}

The following theorem sums up the results of this section which will be used later, with the exception of Lemma~\ref{detuxpos}.

\begin{theo} \label{detavoir}
For $m \in \lbrace 0,1 \rbrace$ and $p \in [1,\infty]$, or for $m \geq 2$ and  $p \in (1,\infty]$, we have:
\begin{equation} \label{detavoirsim}
\Big( \lbrace \left|u(t)\right|_{m,p}^{\alpha} \rbrace \Big)^{1/\alpha} \overset{m,p,\alpha}{\sim} \nu^{-\gamma},\qquad \alpha>0,
\end{equation}
where $\lbrace \cdot \rbrace$ denotes time-averaging over $[T_1,T_2]$. The upper estimates in (\ref{detavoirsim}) hold without time-averaging, uniformly for $t$ separated from $0$. 
Namely, we have
$$
\left|u(t)\right|_{m,p} \overset{m,p}{\lesssim} \max(\nu^{-\gamma},t^{-\gamma}).
$$
On the other hand, the lower estimates hold for all $m \geq 0$, $p \in [1,\infty]$, $\alpha>0$.
\end{theo}

\textbf{Proof.} 
The upper estimates follow from Theorem~\ref{upperwmp}. The lower estimates for $\alpha \geq 2$ follow from Lemma~\ref{detfinalexpbis} by H{\"o}lder's inequality. Finally, for all $m,p$ except $m \geq 1$ and $p=1$ we obtain lower estimates for $\alpha \in (0,2)$ using lower estimates for $\alpha=2$, upper estimates for $\alpha=3$ and H{\"o}lder's inequality. Indeed:
\begin{align} \nonumber
\lbrace \left|u(t)\right|_{m,p}^{\alpha} \rbrace &\geq 
\Big( \lbrace \left|u(t)\right|_{m,p}^{2} \rbrace \Big)^{3-\alpha}
\Big( \lbrace \left|u(t)\right|_{m,p}^{3} \rbrace \Big)^{-(2-\alpha)}
\\ \nonumber
& \gtrsim \nu^{-(6-2\alpha) \gamma} \nu^{(6-3\alpha) \gamma} =\nu^{-\alpha \gamma}.
\end{align}
For $|u|_{m,1},\ m>1$, the lower estimates follow from the ones on $|u|_{m-1,\infty}$.\ $\square$
\medskip
\\ \indent
This theorem yields, for integers $m \geq 1$, the relation
\begin{equation} \label{integers}
\lbrace \Vert u\Vert_m^2 \rbrace \overset{m}{\sim} \nu^{-(2m-1)}.
\end{equation}
By a standard interpolation argument (see (\ref{Sobolevspectr})) the upper bound in (\ref{integers}) also holds for non-integer indices $s >1$. Actually, the same is true for the lower bound, since for any integer $n>s$ we have
\begin{align} \nonumber
\lbrace \Vert u\Vert_s^2 \rbrace &\geq \lbrace \Vert u\Vert_n^2 \rbrace^{n-s+1} \lbrace \Vert u\Vert_{n+1}^2 \rbrace^{-(n-s)} \overset{s}{\gtrsim} \nu^{-(2s-1)}.
\end{align}

\subsection{Estimates for small-scale quantities} \label{detturb}

In this section, we study analogues of quantities which are important for the study of hydrodynamical turbulence. We consider quantities in physical space (structure functions) as well as in Fourier space (energy spectrum). We assume that $\nu \leq \nu_0$. The value of $\nu_0>0$ will be chosen in (\ref{detC1C2}).
\\ \indent
We define the intervals
\begin{equation} \nonumber
J_1=(0,\ C_1 \nu];\ J_2=(C_1 \nu,\ C_2];\ J_3=(C_2,\ 1].
\end{equation}
The nonnegative constants $C_1$ and $C_2$ will be chosen in (\ref{detK}-\ref{detC1C2}) in such a manner that $C_1 \nu_0<C_2<1$, which ensures that the intervals $J_i$ are non-empty and non-intersecting.
\\ \indent
By Theorem \ref{detavoir}, we obtain that $\lbrace |u|^2 \rbrace \sim 1$. On the other hand, by (\ref{detu11}) we get (after integration by parts):
\begin{align} \nonumber
\lbrace |\hat{u}(n)|^2 \rbrace &= (2 \pi n)^{-2} \Big\{ \Big| \int_{S^1}{e^{2 \pi i n x} u_x(x)} \Big|^2 \Big\} 
\\ \label{detintparts}
&\leq (2 \pi n)^{-2} \lbrace |u|^2_{1,1} \rbrace \leq C n^{-2},
\end{align}
and $C_1$ and $C_2$ can be made as small as desired (cf. (\ref{detnu0ineq})). Consequently, the proportion of the sum $\lbrace\sum |\hat{u}(n)|^2 \rbrace$ contained in Fourier modes corresponding to $J_3$ can be made as large as desired. For instance, we may assume that
$$
\Bigg\{ \sum_{|n| < C_2^{-1}}{|\hat{u}(n)|^2} \Bigg\} \geq \frac{99}{100} \Bigg\{ \sum_{n \in \Z}{|\hat{u}(n)|^2} \Bigg\}.
$$
\indent
For $p \geq 0$, we define the structure function of $p$-th order as:
\begin{equation} \nonumber
S_{p}(\ell) = \Big\{ \int_{S^1}{|u(t,x+\ell)-u(t,x)|^p dx} \Big\}.
\end{equation}
The flatness $F(\ell)$, which measures spatial intermittency, is given by
\begin{equation} \label{detflatness}
F(\ell)=S_4(\ell)/S_2^2(\ell).
\end{equation}
Finally, for $k \geq 1$, we define the (layer-averaged) energy spectrum by
\begin{equation} \label{detspectrum}
E(k)=\Bigg\{ \frac{\sum_{|n| \in [M^{-1}k,Mk]}{|\hat{u}(n)|^2}}{\sum_{|n| \in [M^{-1}k,Mk]}{1}} \Bigg\},
\end{equation}
where $M \geq 1$ is a constant which will be specified later (see the proof of Theorem~\ref{detspectrinert}).
\\ \indent
We begin by estimating the functions $S_{p}(\ell)$ from above.

\begin{lemm} \label{detupperdiss}
For $\ell \in [0,1]$,
$$
S_{p}(\ell) \overset{p}{\lesssim} \left\lbrace \begin{aligned} & \ell^{p},\ 0 \leq p \leq 1. \\ & \ell^{p} \nu^{-(p-1)},\ p \geq 1. \end{aligned} \right.
$$
\end{lemm}

\textbf{Proof.} 
We begin by considering the case $p \geq 1$. We have
\begin{align} \nonumber
S_{p}(\ell) &= \Big\{ \int_{S^1}{|u(x+\ell)-u(x)|^p dx} \Big\}
\\ \nonumber
& \leq \Big\{ \Big( \int_{S^1}{|u(x+\ell)-u(x)| dx} \Big) \Big( \max_{x} |u(x+\ell)-u(x)|^{p-1} \Big) \Big\}.
\end{align}
Using the fact that the space average of $u(x+\ell)-u(x)$ vanishes and H{\"o}lder's inequality, we obtain that
\begin{align} \nonumber
S_p(\ell) \leq & \Big\{ \Big(2 \int_{S^1}{(u(x+\ell)-u(x))^{+} dx} \Big)^{p} \Big\}^{1/p} 
\\ \nonumber
& \times \Big\{ \max_{x} |u(x+\ell)-u(x)|^{p} \Big\}^{(p-1)/p}
\\ \label{detdifference}
\leq & C \ell \Big\{ \max_{x} |u(x+\ell)-u(x)|^{p} \Big\}^{(p-1)/p},
\end{align}
where the second inequality follows from Lemma~\ref{detuxpos}. Finally, by Theorem~\ref{detavoir} we get
\begin{align} \nonumber
S_p(\ell) & \leq C \ell \Big\{ ( \ell |u|_{1,\infty} )^{p} \Big\}^{(p-1)/p} \leq
C \ell^{p} \nu^{-(p-1)}.
\end{align}
The case $p<1$ follows immediately from the case $p=1$ since now $S_{p}(\ell) \leq (S_{1}(\ell))^p$, by H{\"o}lder's inequality. $\square$
\medskip
\\ \indent
For $\ell \in J_2 \cup J_3$, we have a better upper bound if $p \geq 1$.

\begin{lemm} \label{detupperinert}
For $\ell \in J_2 \cup J_3$,
$$
S_{p}(\ell) \overset{p}{\lesssim} \left\lbrace \begin{aligned} & \ell^{p},\ 0 \leq p \leq 1. \\ & \ell,\ p \geq 1. \end{aligned} \right.
$$
\end{lemm}

\textbf{Proof.} The calculations are almost the same as in the previous lemma. The only difference is that we use another bound for the right-hand side of (\ref{detdifference}). Namely, by Theorem~\ref{detavoir} we have
\begin{align} \nonumber
S_{p}(\ell) & \leq C \ell \Big\{ \max_{x} |u(x+\ell)-u(x)|^{p} \Big\}^{(p-1)/p}
\\ \nonumber
& \leq C \ell \Big\{ (2 |u|_{\infty} )^{p} \Big\}^{(p-1)/p} \leq C \ell.\ \square
\end{align}

\begin{rmq}
Lemmas~\ref{detupperdiss} and \ref{detupperinert} actually hold even if we drop the time-averaging, since in deriving them we only use upper estimates which hold uniformly for $t \geq T_1$.
\end{rmq}
\indent
To prove the lower estimates for $S_p(\ell)$, we need a lemma. Loosely speaking, this lemma states that there exists a large enough set $L_K \subset [T_1,T_2]$ such that for $t \in L_K$, several Sobolev norms are of the same order as their time averages. Thus, for $t \in L_K$, we can prove the existence of a \enquote{cliff} of height at least $C$ and width at least $C \nu$, using some of the arguments in \cite{AFLV92} which we explained in Subsection~\ref{Burgu}.
\\ \indent
Note that in the following definition, (\ref{detcondi}-\ref{detcondii}) contain lower and upper estimates, while (\ref{detcondiii}) only contains an upper estimate. The inequality $|u(t)|_{\infty} \leq  \max u_x(t)$ in (\ref{detcondi}) always holds, since $u(t)$ has zero mean value and the length of $S^1$ is $1$.

\begin{defi}
For $K>1$, we denote by $L_K$ the set of all $t \in [T_1,T_2]$ such that the assumptions
\begin{align} \label{detcondi}
& K^{-1} \leq |u(t)|_{\infty} \leq  \max u_x(t) \leq K
\\ \label{detcondii}
& K^{-1} \nu^{-1} \leq  |u(t)|_{1,\infty} \leq K \nu^{-1}
\\ \label{detcondiii}
& |u(t)|_{2,\infty} \leq K \nu^{-2}
\end{align}
hold.
\end{defi}

\begin{lemm} \label{dettypical}
There exist constants $C,K_1>0$ such that for $K \geq K_1$, the Lebesgue measure of $L_K$ verifies $\lambda(L_K) \geq C$.
\end{lemm}

\textbf{Proof.}
We begin by noting that if $K \leq K'$, then $L_K \subset L_{K'}$. By Lemma~\ref{detuxpos} and Theorem~\ref{detavoir}, for $K$ large enough the upper estimates in (\ref{detcondi}-\ref{detcondiii}) hold for all $t \geq T_1$. Therefore, if we denote by $B_K$ the set of $t$ such that
$$
\text{\enquote{The lower estimates in (\ref{detcondi}-\ref{detcondii}) hold for a given value of $K$}},
$$
then it suffices to prove the lemma's statement with $B_K$ in place of $L_K$. Now denote by $D_K$ the set of $t$ such that
$$
\text{\enquote{The lower estimate in (\ref{detcondii}) holds for a given value of $K$}}.
$$
By (GN) we have
$$
|u|_{\infty} \geq C |u|_{2,\infty}^{-1} |u|_{1,\infty}^2.
$$
Thus if $D_K$ holds, then $B_{K'}$ holds for $K'$ large enough. Now it remains to show that there exists $C>0$ such that for $K$ large enough, we have the inequality $\lambda(D_K) \geq C$. We clearly have
$$
\lbrace |u|_{1,\infty} \One(|u|_{1,\infty} < K^{-1} \nu^{-1}) \rbrace < K^{-1} \nu^{-1}.
$$
Here, $\One(A)$ denotes the indicator function of an event $A$. On the other hand, by the estimate for $\lbrace |u|_{1,\infty}^2 \rbrace$ in Theorem~\ref{detavoir} we get
\begin{align} \nonumber
\lbrace |u|_{1,\infty} \One( |u|_{1,\infty} > K \nu^{-1}) \rbrace & < K^{-1} \nu \lbrace |u|_{1,\infty}^2 \rbrace \leq C K^{-1} \nu^{-1}
\end{align}
Now denote by $f$ the function
$$
f=|u|_{1,\infty} \One(K_0^{-1} \nu^{-1} \leq |u|_{1,\infty} \leq K_0 \nu^{-1}).
$$
The inequalities above and the lower estimate for $\lbrace |u|_{1,\infty} \rbrace$ in Theorem~\ref{detavoir} imply that
$$
\lbrace f \rbrace \geq (C-K_0^{-1}-C K_0^{-1}) \nu^{-1} \geq C_0 \nu^{-1},
$$
for some suitable constants $C_0$ and $K_0$. Since $f \leq K_0 \nu^{-1}$, we get
$$
\lambda(f \geq C_0 \nu^{-1}/2) \geq C_0 K_0^{-1} (T_2-T_1)/2.
$$
Thus, since $|u|_{1,\infty} \geq f$, we have the inequality
$$
\lambda(|u|_{1,\infty} \geq C_0 \nu^{-1}/2) \geq C_0 K_0^{-1} (T_2-T_1)/2,
$$
which implies the existence of $C,K_1>0$ such that $\lambda(D_{K}) \geq C$ for $K \geq K_1$. $\square$
\smallskip
\\ \indent
Let us denote by $O_K \subset [T_1,T_2]$ the set defined as $L_K$, but with relation (\ref{detcondii}) replaced by
\begin{equation} \label{detcondiibis}
K^{-1} \nu^{-1} \leq -\min u_x \leq K \nu^{-1}.
\end{equation}

\begin{cor} \label{dettypicalcor}
For $K \geq K_1$ and $\nu < K_1^{-2}$, we have $\lambda(O_K) \geq C$. Here, $C,K_1$ are the same as in the formulation of Lemma~\ref{dettypical}.
\end{cor}

\textbf{Proof.} For $K = K_1$ and $\nu < K_1^{-2}$, the estimates (\ref{detcondi}-\ref{detcondii}) tell us that
$$
\max u_x(t) \leq K_1 < K_1^{-1} \nu^{-1} \leq |u_x(t)|_{\infty},\quad t \in L_K.
$$
Thus, in this case we have $O_K=L_K$, which proves the corollary's assertion. Since increasing $K$ while keeping $\nu$ constant increases the measure of $O_K$, it follows for $K \geq K_1$ and $\nu < K_1^{-2}$ we still have $\lambda(O_K) \geq C$. $\square$
\smallskip
\\ \indent
Now we fix
\begin{equation} \label{detK}
K=K_1,
\end{equation}
and choose
\begin{equation} \label{detC1C2}
\nu_0=\frac{1}{6} K^{-2};\ C_1=\frac{1}{4}K^{-2};\ C_2=\frac{1}{20}K^{-4}.
\end{equation}
In particular, we have $0<C_1 \nu_0 <C_2<1$: thus the intervals $J_i$ are non-empty and non-intersecting for all $\nu \in (0,\nu_0]$. Everywhere below the constants depend on $K$.
\\ \indent
Actually, we can choose any values of $C_1$, $C_2$ and $\nu_0$, provided
\begin{equation} \label{detnu0ineq}
C_1 \leq \frac{1}{4}K^{-2};\quad 5 K^2 \leq \frac{C_1}{C_2}<\frac{1}{\nu_0}.
\end{equation}

\begin{lemm} \label{detlowerdiss}
For $\ell \in J_1$,
$$
S_{p}(\ell) \overset{p}{\gtrsim} \left\lbrace \begin{aligned} &\ell^{p},\ 0 \leq p \leq 1. \\ & \ell^{p} \nu^{-(p-1)},\ p \geq 1. \end{aligned} \right.
$$
\end{lemm}

\textbf{Proof.}
By Corollary~\ref{dettypicalcor}, it suffices to prove that the inequalities hold uniformly in $t$ for $t \in O_K$, with $S_{p}(\ell)$ replaced by
\begin{equation} \nonumber
\int_{S^1}{|u(x+\ell)-u(x)|^p dx}.
\end{equation}
Till the end of this proof, we assume that $t \in O_K$.
\\ \indent
Denote by $z$ the leftmost point on $S^1$ (considered as $[0,1)$) such that $ u'(z) \leq - K^{-1} \nu^{-1}$. Since $|u|_{2,\infty} \leq K \nu^{-2}$, we have
\begin{equation} \label{detuxsmall}
u'(y) \leq -\frac{1}{2} K^{-1}  \nu^{-1},\quad y \in [z-\frac{1}{2} K^{-2} \nu,z+\frac{1}{2} K^{-2} \nu].
\end{equation}
In other words, the interval
$$
[z-\frac{1}{2} K^{-2} \nu,z+\frac{1}{2} K^{-2} \nu]
$$
corresponds to (a part of) a cliff.
\\ \indent
\textbf{Case $\mathbf{p \geq 1}$.} Since $\ell \leq C_1 \nu=\frac{1}{4}K^{-2} \nu$, by H{\"o}lder's inequality we get
\begin{align} \nonumber
\int_{S^1}&{|u(x+\ell)-u(x)|^p dx} \geq \int_{z-\frac{1}{4} K^{-2} \nu}^{z+\frac{1}{4} K^{-2} \nu}{|u(x+\ell)-u(x)|^p dx}
\\ \nonumber
&\geq (K^{-2} \nu/2)^{1-p} \Big( \int_{z-\frac{1}{4} K^{-2} \nu}^{z+\frac{1}{4} K^{-2} \nu}{|u(x+\ell)-u(x)| dx} \Big)^p
\\ \nonumber
&= C(p) \nu^{1-p} \Big( \int_{z-\frac{1}{4} K^{-2} \nu}^{z+\frac{1}{4} K^{-2} \nu}{ \Big(\int_{x}^{x+\ell}{- u'(y) dy } \Big) dx} \Big)^p
\\ \nonumber
&\geq C(p) \nu^{1-p} \Big( \int_{z-\frac{1}{4} K^{-2} \nu}^{z+\frac{1}{4} K^{-2} \nu}{ \frac{1}{2} \ell K^{-1} \nu^{-1} \ dx} \Big)^p = C(p) \nu^{1-p} \ell^p.
\end{align}
\\ \indent
\textbf{Case $\mathbf{p < 1}$.} By H{\"o}lder's inequality we obtain that
\begin{align} \nonumber
&\int_{S^1}{|u(x+\ell)-u(x)|^p dx} \geq \int_{S^1}{\Big((u(x+\ell)-u(x))^+\Big)^p dx}
\\ \nonumber
&\geq \Big( \int_{S^1}{\Big((u(x+\ell)-u(x))^+ \Big)^2 dx} \Big)^{p-1} \Big( \int_{S^1}{(u(x+\ell)-u(x))^+ dx} \Big)^{2-p}.
\end{align}
Using the upper estimate in (\ref{detcondi}) we get
\begin{align} \nonumber
&\int_{S^1}{|u(x+\ell)-u(x)|^p dx}
\\ \nonumber
&\geq \Big( \int_{S^1}{\ell^2 K^2 dx} \Big)
^{p-1} \Big( \int_{S^1}{(u(x+\ell)-u(x))^+ dx} \Big)^{2-p}.
\end{align}
Since $\int_{S^1}{(u(\cdot+\ell)-u(\cdot))}=0$, we obtain that
\begin{align} \nonumber
&\int_{S^1}{|u(x+\ell)-u(x)|^p dx}
\\ \nonumber
&\geq C(p) \ell^{2 (p-1)} \Big(\frac{1}{2} \int_{S^1}{|u(x+\ell)-u(x)| dx} \Big)^{2-p} \geq C(p) \ell^p.
\end{align}
The last inequality follows from the case $p=1$.\ $\square$
\smallskip
\\ \indent
The proof of the following lemma uses an argument from \cite{AFLV92}, which becomes quantitative if we restrict ourselves to the set $O_K$.

\begin{lemm} \label{detlowerinert}
For $m \geq 0$ and $\ell \in J_2$,
$$
S_{p}(\ell) \overset{p}{\gtrsim} \left\lbrace \begin{aligned} & \ell^{p},\ 0 \leq p \leq 1. \\ & \ell,\ p \geq 1. \end{aligned} \right.
$$
\end{lemm}

\textbf{Proof.} In the same way as above, it suffices to prove that the inequalities hold uniformly in $t$ for $t \in O_K$, with $S_{p}(\ell)$ replaced by
\begin{equation} \nonumber
\int_{S^1}{|u(x+\ell)-u(x)|^p dx},
\end{equation}
and we can restrict ourselves to the case $p \geq 1$. Again, till the end of this proof, we assume that $t \in O_K$.
\\ \indent
Define $z$ as in the proof of Lemma~\ref{detlowerdiss}. We have
\begin{align} \nonumber
\int_{S^1}&{|u(x+\ell)-u(x)|^p dx} \geq
\\ \nonumber
&\int_{z-\frac{1}{2}\ell}^{z}
{ \Big| \underbrace{\int_{x}^{x+\ell}{u'^-(y)dy}}_{"cliffs"} - \underbrace{\int_{x}^{x+\ell}{u'^+(y)dy}}_{"ramps"} \Big|^p dx}.
\end{align}
Since $\ell \geq C_1 \nu=\frac{1}{4} K^{-2} \nu$, by (\ref{detuxsmall}) for $x \in [z-\frac{1}{2}\ell,z]$ we get
\begin{align} \nonumber
\int_{x}^{x+\ell}{u'^-(y)dy} &\geq \int_{z}^{z+\frac{1}{8} K^{-2} \nu}{u'^-(y) dy} \geq
\frac{1}{16} K^{-3}.
\\ \nonumber
&.
\end{align}
On the other hand, since $\ell \leq C_2$, by (\ref{detcondi}) and (\ref{detC1C2}) we have
$$
\int_{x}^{x+\ell}{u'^+(y)dy} \leq C_2 K = \frac{1}{20} K^{-3}.
$$
Thus,
\begin{align} \nonumber
& \int_{S^1}{|u(x+\ell)-u(x)|^p dx} \geq \frac{1}{2} \ell \Bigg( \Big(\frac{1}{16}-\frac{1}{20}\Big) K^{-3} \Bigg)^p \geq C(p) \ell.\ \square
\end{align}
\medskip \\ \indent
Summing up the results above we obtain the following theorem.

\begin{theo} \label{detavoir2}
For $\ell \in J_1$,
$$
S_{p}(\ell) \overset{p}{\sim} \left\lbrace \begin{aligned} & \ell^{p},\ 0 \leq p \leq 1. \\ & \ell^{p} \nu^{-(p-1)},\ p \geq 1. \end{aligned} \right.
$$
On the other hand, for $\ell \in J_2$,
$$
S_{p}(\ell) \overset{p}{\sim} \left\lbrace \begin{aligned} & \ell^{p},\ 0 \leq p \leq 1. \\ & \ell,\ p \geq 1. \end{aligned} \right.
$$
\end{theo}

The following result follows immediately from the definition (\ref{detflatness}).

\begin{cor} \label{detflatnesscor}
For $\ell \in J_2$, the flatness satisfies $F(\ell) \sim \ell^{-1}$.
\end{cor}

By Theorem~\ref{detavoir}, for $m \geq 1$ we have
$$
\lbrace |\hat{u}(k)|^2 \rbrace \leq (2 \pi k)^{-2m} \lbrace \Vert u \Vert_m^2 \rbrace \overset{m}{\sim} (k \nu)^{-2m} \nu.
$$
Thus, for $|k| \succeq \nu^{-1}$, $\lbrace|\hat{u}(k)|^2\rbrace$ decreases super-algebraically.
\smallskip
\\ \indent
Now we want to estimate the $H^s$ norms of $u$ for $s \in (0,1)$.

\begin{lemm} \label{detH12}
We have
$$
\lbrace \Vert u\Vert_{1/2}^2 \rbrace \sim |\log \nu|.
$$
\end{lemm}

\textbf{Proof.}
By (\ref{Sobolevfrac}) we have
\begin{align} \nonumber
\left\|u\right\|_{1/2} \sim \Bigg( \int_{S^1} \Big(\int_0^1 {\frac{|u(x+\ell)-u(x)|^2}{\ell^{2}} d \ell} \Big) dx \Bigg)^{1/2}.
\end{align}
Consequently, by Fubini's theorem,
\begin{align} \nonumber
\\ \nonumber
\lbrace \left\|u\right\|^2_{1/2} \rbrace &\sim \int_0^1 \frac{1}{\ell^2} \Big\{ \int_{S^1}{|u(x+\ell)-u(x)|^2 dx} \Big\} d \ell 
\\ \nonumber
&= \int_{0}^{1}{\frac{S_2(\ell)}{\ell^2} d\ell}
=\int_{J_1}{\frac{S_2(\ell)}{\ell^2} d\ell}+\int_{J_2}{\frac{S_2(\ell)}{\ell^2} d\ell}+\int_{J_3}{\frac{S_2(\ell)}{\ell^2} d\ell}.
\end{align}
By Theorem~\ref{detavoir2} we get
$$
\int_{J_1}{\frac{S_2(\ell)}{\ell^2} d\ell} \sim \int_{0}^{C_1 \nu}{\frac{\ell^2 \nu^{-1}}{\ell^2} d\ell} \sim 1
$$
and
$$
\int_{J_2}{\frac{S_2(\ell)}{\ell^2} d\ell} \sim \int_{C_1 \nu}^{C_2}{\frac{\ell}{\ell^2} d\ell} \sim |\log  \nu|,
$$
respectively. Finally, by Lemma~\ref{detupperinert} we get
$$
\int_{J_3}{\frac{S_2(\ell)}{\ell^2} d\ell} \leq C C_2^{-2} \leq C.
$$
Thus,
$$
\lbrace \left\|u\right\|^2_{1/2} \rbrace \sim |\log \nu|.\ \qed
$$
\smallskip \\ \indent
The proof of the following result follows the same lines.

\begin{lemm} \label{H01}
For $s \in (0,1/2)$,
$$
\lbrace \Vert u\Vert_s^2 \rbrace \overset{s}{\sim} 1.
$$
On the other hand, for $s \in (1/2,1)$,
$$
\lbrace \Vert u\Vert_s^2 \rbrace \overset{s}{\sim} \nu^{-(2s-1)}.
$$
\end{lemm}

The results above tell us that $\lbrace|\hat{u}(k)|^2\rbrace$ decreases very fast for $|k| \gtrsim \nu^{-1}$ and that for $s \geq 0$ the sums $\sum{|k|^{2s} \lbrace |\hat{u}(k)|^2}\rbrace$ have exactly the same behaviour as the partial sums $\sum_{|k| \leq \nu^{-1}}{|k|^{2s} |k|^{-2}}$ in the limit $\nu \rightarrow 0^+$. Therefore we can conjecture that for $|k| \lesssim \nu^{-1}$, we have  $\lbrace|\hat{u}(k)|^2 \rbrace \sim |k|^{-2}$.
\\ \indent
A result of this type actually holds (after layer-averaging), as long as $|k|$ is not too small. To prove it, we use a version of the Wiener--Khinchin theorem, stating that for any function $v \in L_2$ one has
\begin{equation} \label{detspectrinertaux}
|v(\cdot+y)-v(\cdot)|^2=4\sum_{n \in \Z}{ \sin^2 (\pi ny) |\hat{v}(n)|^2}.
\end{equation}

\begin{theo} \label{detspectrinert}
For $k$ such that $k^{-1} \in J_2$, we have $E(k) \sim k^{-2}$.
\end{theo}

\textbf{Proof.}
We recall that by definition (\ref{detspectrum}),
$$
E(k) = \Bigg\{ \frac{\sum_{|n| \in [M^{-1}k,Mk]}{|\hat{u}(n)|^2}}{\sum_{|n| \in [M^{-1}k,Mk]}{1}} \Bigg\}.
$$
Therefore proving the assertion of the theorem is the same as proving that
\begin{equation} \label{detspectrinertequiv}
\sum_{|n| \in [M^{-1}k,Mk]}{ n^2 \lbrace |\hat{u}(n)|^2 \rbrace } \sim k.
\end{equation}
From now on, we will indicate explicitly the dependence on $M$. The upper estimate holds without averaging over $n$ such that
$$
|n| \in [M^{-1} k, Mk].
$$
Indeed,  by (\ref{detintparts}) we know that
\begin{equation} \nonumber
\lbrace |\hat{u}(n)|^2 \rbrace \leq C n^{-2}.
\end{equation}
Also, this inequality implies that
\begin{equation} \label{detspectrinertupper1}
\sum_{|n| < M^{-1} k}{ n^2 \lbrace |\hat{u}(n)|^2 \rbrace } \leq C M^{-1} k
\end{equation}
and
\begin{equation} \label{detspectrinertupper2}
\sum_{|n| > M k}{\lbrace |\hat{u}(n)|^2 \rbrace } \leq C M^{-1} k^{-1}.
\end{equation}
To prove the lower bound we note that
\begin{align} \nonumber
\sum_{|n| \leq M k}{ n^2 \lbrace |\hat{u}(n)|^2 \rbrace } & \geq \frac{k^2}{\pi^2} \sum_{|n| \leq M k}{ \sin^2 (\pi nk^{-1}) \lbrace |\hat{u}(n)|^2 \rbrace }
\\ \nonumber
&\geq \frac{k^2}{\pi^2} \Big( \sum_{n \in \Z}{ \sin^2 (\pi nk^{-1}) \lbrace |\hat{u}(n)|^2 \rbrace } - \sum_{|n| > M k}{\lbrace |\hat{u}(n)|^2 \rbrace } \Big).
\end{align}
Using (\ref{detspectrinertaux}) and (\ref{detspectrinertupper2}) we get
\begin{align} \nonumber
\sum_{|n| \leq M k}{ n^2 \lbrace |\hat{u}(n)|^2 \rbrace } &\geq  \frac{k^2}{4 \pi^2} \Big( \lbrace |u(\cdot+k^{-1})-u(\cdot)|^2 \rbrace - C M^{-1} k^{-1} \Big)
\\ \nonumber
&\geq \frac{k^2}{4 \pi^2} \Big( S_2(k^{-1})-C M^{-1} k^{-1} \Big).
\end{align}
Finally, using Theorem~\ref{detavoir2} we obtain that
\begin{equation} \nonumber
\sum_{|n| \leq M k}{ n^2 \lbrace |\hat{u}(n)|^2 \rbrace }  \geq (C-C M^{-1}) k.
\end{equation}
Now we use (\ref{detspectrinertupper1}) and we choose $M \geq 1$ large enough to obtain (\ref{detspectrinertequiv}). $\square$

\section{The randomly forced Burgers equation} \label{rand}

\subsection{Foreword}

The results stated in this section have been obtained in \cite{BorW} for the white-forced equation. For the simpler case of the kick force, estimates for Sobolev norms have been obtained in \cite{BorK}. Since those estimates are used as a \enquote{black box} when studying small-scale quantities, generalisation of the small-scale estimates in \cite{BorW} to the case of a kick force is immediate. Thus, in this section, we only consider the white-forced equation (\ref{whiteBurgers}).
\\ \indent
Existence and uniqueness of smooth solutions to (\ref{whiteBurgers}) is proved by the "mild solution" technique  (cf. \cite[Chapter 14]{DZ96}). For the kicked equation, existence and uniqueness of solutions follows from the corresponding fact for the unforced equation.
\\ \indent
Some proofs in \cite{BorW} are similar to the proofs in the unforced case. We will only give here the proofs of Theorem~\ref{whiteuxpos} and Lemma~\ref{whitefinitetime}, as well as some comments on the proofs of small-scale results.
\\ \indent
The major difference between the unforced and the white-forced generalised Burgers equation is the energetic picture. In the first case, we have a dissipative system: the $L_2$ norm is decreasing in time. Consequently, the regime where energy dissipates fast enough (which yields a time-averaged lower bound on the Sobolev norms) is transient and depends on the initial condition. On the contrary, in the second case, after a time needed either to dissipate energy if $u_0$ is large or to supply energy if $u_0$ is small, we are in a \textit{quasi-stationary regime}, in the sense that in average on large enough time intervals, we have an approximate balance between the dissipation rate $-\nu \E \Vert u \Vert_1^2$ and the constant energy supply rate $I_0$.
\\ \indent
For simplicity, in the white-forced case we assume that the initial condition $u_0$ is deterministic. However, we can easily generalise all results to the case of a random initial condition independent of $w(t), t \geq 0$. Indeed, in that case for any measurable functional $\Phi(u(\cdot))$ we have
$$
\E \Phi(u(\cdot))=\int{\E \Big(\Phi(u(\cdot))|u(0)=u_0 \Big) d \mu(u_0)},
$$
where $\mu(u_0)$ is the law of $u_0$, and all our estimates hold uniformly in $u_0$.
\\ \indent
Moreover, for $\tau \geq 0$ and $u_0$ independent of $w(t)-w(\tau), t \geq \tau$,  the Markov property yields:
$$
\E \Phi(u(\cdot))=\int{\E \Big(\Phi(u(\tau+\cdot))|u(\tau)=u_0 \Big) d \mu(u_0)}.
$$
Consequently, all estimates which hold for time $t$ or a time interval $[t,t+T]$ actually hold for time $t+\tau$ or a time interval $[t+\tau,t+\tau+T]$, uniformly in $\tau \geq 0$.
\\ \indent
The remarks above still hold for the kick-forced equation. However, the constant energy supply rate (and continuous time-invariance of the forcing) are replaced by constant energy supply at the discrete moments $i \in \N$ (and discrete time-invariance of the forcing).

\subsection{Estimates for Sobolev norms}

The following theorem is proved using a stochastic version of the Kruzhkov maximum principle (cf. \cite{Kru64}). In all results in this section, quantities estimated for fixed $\omega$, such as $\max_{s \in [t,t+1],\ x \in S^1}{u_x}$ or maxima in time of Sobolev norms, can be replaced by their suprema over all smooth initial conditions. For instance, the quantity
$$
\max_{s \in [t,t+1]} |u(s)|_{m,p}
$$
can be replaced by
$$
\sup_{u_0 \in C^{\infty}} \max_{s \in [t,t+1]} |u(s)|_{m,p}.
$$
For the lower estimates, this fact is obvious. For the upper ones, the reason is that these quantities admit upper bounds of the form
$$
(1+\max_{s \in [t-\tau,t+\tau]}{\Vert w(s) \Vert_m})^{\alpha(m)} \nu^{-\beta(m)}.
$$

\begin{theo} \label{whiteuxpos}
Denote by $X_t$ the random variable
$$
X_t=\max_{s \in [t,t+1]} \max_{x \in S^1} u_x(s,x).
$$
For every $k \geq 1$, we have
$$
\E \ X_t^{k} \overset{k}{\lesssim} 1,\quad t \geq 1.
$$
\end{theo}

\textbf{Proof.}
We take $t=1$, denoting $X_t$ by $X$.
\\ \indent
Consider (\ref{whiteBurgers}) on the time interval $[0,2]$. Putting $v=u-w$ and differentiating once in space,  we get
\begin{equation} \label{whitemaxvx}
\frac{\partial v_x}{\partial t} + f''(u)(v_x+w_x)^2+f'(u) (v_{x}+w_{x})_x=\nu (v_{x}+w_{x})_{xx}.
\end{equation}
Consider $\tilde{v}(t,x)=tv_x(t,x)$ and multiply (\ref{whitemaxvx}) by $t^2$. For $t>0$, $\tilde{v}$ verifies
\begin{align} \nonumber
&t\tilde{v}_t -\tilde{v} + f''(u) (\tilde{v}+tw_x)^2 + tf'(u) \tilde{v}_x + t^2 f'(u) w_{xx}
\\ \label{whitemaxtildev}
&= \nu t \tilde{v}_{xx}+\nu t^2 w_{xxx}.
\end{align}
Now observe that if the zero mean function $\tilde{v}$ does not vanish identically on the domain $S=\left[0,2\right] \times S^1$, then it attains its positive maximum $N$ on $S$ at a point $(t_1,x_1)$ such that $t_1>0$. At $(t_1,x_1)$ we have $\tilde{v}_t \geq 0$, $\tilde{v}_x=0$ and $\tilde{v}_{xx} \leq 0$.
By (\ref{whitemaxtildev}), at $(t_1,x_1)$ we have the inequality
\begin{equation} \label{whitemaxpoint}
f''(u) (\tilde{v}+tw_x)^2 \leq \tilde{v}- t^2 f'(u) w_{xx}+ \nu t^2 w_{xxx}.
\end{equation}
Denote by $A$ the random variable
$$
A=\max_{t \in [0,2]} |w(t)|_{3,\infty}.
$$
Since for every $t$, $tv(t)$ is the zero space average primitive of $\tilde{v}(t)$ on $S^1$, we get
\begin{align} \nonumber
\max_{t \in [0,2],\ x \in S^1}{|tu|} &\leq \max_{t \in [0,2],\ x \in S^1}{(|tv|+|tw|)} 
\\ \nonumber
&\leq N+2\max_{t \in [0,2]} |w(t)|_{\infty} \leq N+2A.
\end{align}
Now denote by $\delta$ the quantity
$$
\delta=2-h(1).
$$
By (\ref{poly}), $\delta>0$. We obtain that
\begin{align} \nonumber
\max_{t \in [0,2],\ x \in S^1} |t^2 f'(u) w_{xx}| & \leq A \max_{t \in [0,2],\ x \in S^1} {t^{\delta} |t^{2-\delta} f'(u)|}
\\ \nonumber
&\leq C A \max_{t \in [0,2],\ x \in S^1} {t^{\delta} (|tu|+t)^{2-\delta}}
\\ \nonumber
&\leq C A (N+2A+2)^{2-\delta}.
\end{align}
From now on, we assume that $N \geq 2A$. Since $\nu \in (0,1]$ and $f'' \geq \sigma$, the relation (\ref{whitemaxpoint}) yields
$$
\sigma (N-2A)^2 \leq N+C A (N+2A+2)^{2-\delta}+4A.
$$
Thus we have proved that if $N \geq 2A$, then $N \leq C(A+1)^{1/\delta}$. Since by (\ref{moments}), all moments of $A$ are finite, all moments of $N$ are also finite. By definition of $\tilde{v}$ and $S$, the same is true for $X$. This proves the theorem's assertion. $\qed$

\begin{cor} \label{whiteW11}
For $k \geq 1$,
$$
\E \max_{s \in [t,t+1]} \left|u(s)\right|^k_{1,1} \overset{k}{\lesssim} 1,\quad t \geq 1.
$$
\end{cor}

\begin{cor} \label{whiteLpupper}
For $k \geq 1$,
$$
\E \max_{s \in [t,t+1]} \left|u(s)\right|^k_{p} \overset{k}{\lesssim} 1,\quad p \in [1,\infty],\ t \geq 1.
$$
\end{cor}

\begin{lemm} \label{whiteuppermlemm}
For $m \geq 1$,
$$
\E  \max_{s \in [t,t+1]} \left\|u(s)\right\|^{2}_m \overset{m}{\lesssim} \nu^{-(2m-1)},\quad t \geq 2.
$$
\end{lemm}

\begin{theo} \label{whiteupperwmp}
For $m \in \lbrace 0,1 \rbrace$ and $p \in [1,\infty]$, or for $m \geq 2$ and $p \in (1,\infty]$,
$$
\Big( \E  \max_{s \in [t,t+1]} \left|u(s)\right|^{\alpha}_{m,p} \Big)^{1/\alpha} \overset{m,p,\alpha}{\lesssim} \nu^{-\gamma},\quad \alpha>0,\ t \geq 2.
$$
\end{theo}

\begin{lemm} \label{whitefinitetime}
There exists a constant $T_0>0$ such that we have
$$
\Big( \frac{1}{T} \int_{t}^{t+T}{\ \E \left\|u(s)\right\|_1^2} \Big)^{1/2} \gtrsim \nu^{-1/2},\qquad t \geq 1,\ T \geq T_0.
$$
\end{lemm}

\textbf{Proof.} For $T>0$, by (\ref{Itoexp0diff}) we get
\begin{align} \nonumber
\E \left|u(t+T)\right|^2 &\geq \E ( \left|u(t+T)\right|^2-\left|u(t)\right|^2 )
= TI_0 - 2 \nu \int_{t}^{t+T}{\E \left\|u(s)\right\|_{1}^2}.
\end{align}
On the other hand, by Corollary~\ref{whiteLpupper} there exists a constant $C'>0$ such that
$\E \left|u(t+T)\right|^2 \leq C'$. Consequently, for $T \geq T_0:=(C'+1)/I_0$,
$$
\frac{1}{T} \int_{t}^{t+T}{\E \left\|u(s)\right\|_1^2} \geq \frac{TI_0-C'}{2T} \nu^{-1} \geq \frac{I_0}{2(C'+1)} \nu^{-1},
$$
which proves the lemma's assertion.\ $\qed$

\begin{theo} \label{whiteavoir}
For $m \in \lbrace 0,1 \rbrace$ and $p \in [1,\infty]$, or for $m \geq 2$ and $p \in (1,\infty]$, we have
\begin{align} \nonumber
&\Big( \frac{1}{T} \int_{t}^{t+T}{\E \left|u(s)\right|_{m,p}^{\alpha}} \Big)^{1/\alpha} \overset{m,p,\alpha}{\sim} \nu^{-\gamma},\ \alpha>0,\ 
\\ \label{whiteasymp}
& t \geq T_1=T_0+2,\ T \geq T_0.
\end{align}
Moreover, the upper estimates hold with time-averaging replaced by maximising over $[t,t+1]$, i.e.
\begin{equation} \label{whitemaxim}
\Big( \E \max_{s \in [t,t+1]}{\left|u(s)\right|_{m,p}^{\alpha}}\Big)^{1/\alpha} \overset{m,p,\alpha}{\lesssim} \nu^{-\gamma},\quad \alpha>0,\ t \geq 2.
\end{equation}
On the other hand, the lower estimates hold for all $m \geq 0$ and $p \in [1,\infty]$. The asymptotics (\ref{whiteasymp}) hold without time-averaging if $m$ and $p$ are such that $\gamma(m,p)=0$. Namely, in this case,
\begin{equation} \label{whitegamma0}
\Big( \E \left|u(t)\right|_{m,p}^{\alpha}\Big)^{1/\alpha} \overset{m,p,\alpha}{\sim} 1,\quad \alpha>0,\ t \geq T_1.
\end{equation}
Finally, note that all these estimates hold if we replace Sobolev norms with their suprema over all smooth initial conditions.
\end{theo}

\subsection{Estimates for small-scale quantities} \label{whiteturb}

Consider an observable $A$, i.e. a real-valued functional on a Sobolev space $H^m$, which we evaluate on the solutions $u^{\omega}(s)$. We denote by $\lbrace A \rbrace$ the average of $A(u^{\omega}(s))$ in ensemble and in time over $[t,t+T_0]$:
$$
\lbrace A \rbrace=\frac{1}{T_0}\ \int_{t}^{t+T_0}{\E A(u^{\omega}(s)) ds},\ t \geq T_1.
$$
The constant $T_1$ is the same as in Theorem~\ref{whiteavoir}. In this section, we assume that $\nu \leq \nu_0$, where $\nu_0$ is a nonnegative constant. The definitions and the choices for $\nu_0$, the ranges and the small-scale quantities are word-to-word the same as in the unforced case, up to the changes in the meaning of the brackets $\lbrace \cdot \rbrace$.

\begin{lemm} \label{whiteupperdiss}
For $\alpha \geq 0$ and $\ell \in [0,1]$,
$$
S_{p,\alpha}(\ell) \overset{p,\alpha}{\lesssim} \left\lbrace \begin{aligned} & \ell^{\alpha p},\ 0 \leq p \leq 1. \\ & \ell^{\alpha p} \nu^{-\alpha(p-1)},\ p \geq 1. \end{aligned} \right.
$$
\end{lemm}

\begin{lemm} \label{whiteupperinert}
For $\alpha \geq 0$ and $\ell \in J_2 \cup J_3$,
$$
S_{p,\alpha}(\ell) \overset{p,\alpha}{\lesssim} \left\lbrace \begin{aligned} & \ell^{\alpha p},\ 0 \leq p \leq 1. \\ & \ell^{\alpha},\ p \geq 1. \end{aligned} \right.
$$
\end{lemm}

The following lemma states that with a probability which is not too small, during a period of time which is not too small, several Sobolev norms are of the same order as their expected values.

\begin{defi}
For a given solution $u(s)=u^{\omega}(s)$ and $K>1$, we denote by $L_K$ the set of all $(s,\omega) \in [t,t+T_0] \times \Omega$ such that
\begin{align} \label{whitecondi}
&K^{-1} \leq |u(s)|_{\infty} \leq \max u_x(s) \leq K
\\ \label{whitecondii}
&K^{-1} \nu^{-1} \leq |u(s)|_{1,\infty} \leq K \nu^{-1}
\\ \label{whitecondiii}
&|u(s)|_{2,\infty} \leq K \nu^{-2}.
\end{align}
\end{defi}

\begin{lemm} \label{whitetypical}
There exist constants $\tilde{C},K_1>0$ such that for all $K \geq K_1$, $\rho(L_K) \geq \tilde{C}$. Here, $\rho$ denotes the product measure of the Lebesgue measure and $\Pe$ on $[t,t+T_0] \times \Omega$.
\end{lemm}

\textbf{Proof.}
The proof is almost the same as in the deterministic case. One difference is that now we average in time and in probability instead of only averaging in time. The other difference is that the upper estimates now hold with probability tending to $1$ as $K \rightarrow +\infty$, and not with probability $1$ for $K$ large enough.\ $\qed$

\begin{defi}
For a given solution $u(s)=u^{\omega}(s)$ and $K>1$, we denote by $O_K$ the set of all $(s,\omega) \in [t,t+T_0] \times \Omega$ such that the conditions (\ref{whitecondi}), (\ref{whitecondiii}) and
\begin{equation} \label{whitecondii'}
K^{-1} \nu^{-1} \leq -\min u_x \leq K \nu^{-1}
\end{equation}
hold.
\end{defi}

\begin{cor} \label{whitetypicalcor}
If $K \geq K_1$ and $\nu < K_1^{-2}$, then $\rho(O_K) \geq \tilde{C}$. Here, $\tilde{C}$ and $K_1$ are the same as in the statement of Lemma~\ref{whitetypical}.
\end{cor}

\begin{theo} \label{whiteavoir2}
For $\alpha \geq 0$ and $\ell \in J_1$,
$$
S_{p,\alpha}(\ell) \overset{p,\alpha}{\sim} \left\lbrace \begin{aligned} & \ell^{\alpha p},\ 0 \leq p \leq 1. \\ & \ell^{\alpha p} \nu^{-\alpha (p-1)},\ p \geq 1. \end{aligned} \right.
$$
On the other hand, for $\alpha \geq 0$ and $\ell \in J_2$,
$$
S_{p,\alpha} (\ell) \overset{p,\alpha}{\sim} \left\lbrace \begin{aligned} & \ell^{\alpha p},\ 0 \leq p \leq 1. \\ & \ell^{\alpha},\ p \geq 1. \end{aligned} \right.
$$
\end{theo}

\begin{cor} \label{whiteflatnesscor}
For $\ell \in J_2$, the flatness satisfies $F(\ell) \sim \ell^{-1}$.
\end{cor}

\begin{lemm} \label{whiteH01}
We have
\begin{align} \nonumber
&\lbrace \Vert u\Vert_s^2 \rbrace \overset{s}{\sim} 1,\quad s \in (0,1/2).
\\ \nonumber
& \lbrace \Vert u\Vert_{1/2}^2 \rbrace \sim |\log \nu|.
\\ \nonumber
& \lbrace \Vert u\Vert_s^2 \rbrace \overset{s}{\sim} \nu^{-(2s-1)},\quad s \in (1/2,1).
\end{align}
\end{lemm}

\begin{theo} \label{whitespectrinert}
If $M$ in the definition of $E(k)$ is large enough, then for every $k$ such that $k^{-1} \in J_2$, we have $E(k) \sim k^{-2}$.
Moreover, we have
$$
\Bigg\{ \Bigg( \frac{\sum_{|n| \in [M^{-1}k,Mk]}{|\hat{u}(n)|^2}}{\sum_{|n| \in [M^{-1}k,Mk]}{1}} \Bigg)^{\alpha} \Bigg\} \overset{\alpha}{\sim} k^{-2\alpha},\quad \alpha>0.
$$
\end{theo}

\section{Stationary measure and related issues} \label{stat}

The results in this section are proved in \cite{BorW} for the equation with white forcing. Up to some changes, they can be generalised to the kick force case. For more details, see \cite{BorPhD}; see also \cite{KuSh12}, where a random forcing is introduced in a similar setup.

\begin{theo} \label{contract}
Consider two solutions $u$, $\overline{u}$ of (\ref{whiteBurgers}), corresponding to the same random force but different initial conditions in $C^{\infty}$. For all $t \geq 0$, we have
\begin{equation} \nonumber
|u(t)-\overline{u}(t)|_{1} \leq |u(0)-\overline{u}(0)|_{1}.
\end{equation}
\end{theo}
\indent
Since $C^{\infty}$ is dense in $L_1$, Theorem~\ref{contract} allows us to define solutions of (\ref{whiteBurgers}) for any initial condition in $L_1$. In the same way as in the case of a smooth initial condition, we can prove that those solutions make a time-continuous Markov process, and then we can define the corresponding semigroup $S_t^{*}$ acting on Borel measures on $L_1$. For a more detailed account on the well-posedness in a similar setting, see \cite{KuSh12}.
\\ \indent
A \textit{stationary measure} is a Borel probability measure on $L_1$ invariant by $S_t^{*}$ for every $t$.  A \textit{stationary solution} of (\ref{whiteBurgers}) is a random process $v$ defined for $(t,\omega) \in [0,+\infty) \times \Omega$ and valued in $L_1$, which verifies (\ref{whiteBurgers}), such that the distribution of $v(t,\cdot)$ does not depend on $t$. This distribution is automatically a stationary measure.
\\ \indent
It remains to show existence and uniqueness of a stationary measure, which implies existence and uniqueness (in the sense of distribution) of a stationary solution. Moreover, we obtain an additional bound for the rate of convergence to the stationary measure in an appropriate distance. This bound holds independently from the viscosity or from the initial condition.

\begin{defi}
Fix $p \in [1,\infty)$. For a continuous real-valued function $g$ on $L_p$, we define its Lipschitz norm as
$$
|g|_{L}:=\sup_{L_p}{|g|}+|g|_{Lip},
$$
where $|g|_{Lip}$ is the Lipschitz constant of $g$. The set of continous functions with finite Lipschitz norm will be denoted by $L=L(L_p)$. The choice of $p$ will always be clear from the context.
\end{defi}

\begin{defi}
For two Borel probability measures $\mu_1,\mu_2$ on $L_p$, we denote by $\Vert \mu_1-\mu_2 \Vert^*_L$ the Lipschitz-dual distance:
$$
\Vert \mu_1-\mu_2 \Vert^*_L:=\sup_{g \in L,\ |g|_{L} \leq 1}{\Big| \int_{S^1}{g d \mu_1}-\int_{S^1}{g d \mu_2} \Big|}.
$$
\end{defi}

Since we have $u_0$-uniform upper estimates, existence of a stationary measure for the generalised Burgers equation is proved using the Bogolyubov-Krylov argument (see \cite{KuSh12}).
\\ \indent
Now we state the main result of this section. It immediately implies uniqueness of a stationary measure $\mu$ for the equation (\ref{whiteBurgers}).

\begin{theo} \label{algCV}
There exists a positive constant $C'$ such that for $t \geq 0$, we have
\begin{equation} \label{algCVformula}
\Vert S_t^{*}\mu_1-S_t^{*}\mu_2 \Vert^*_L \leq C't^{-1/13},\qquad t \geq 1,
\end{equation}
for any probability measures $\mu_1$, $\mu_2$ on $L_1$.
\end{theo}

\begin{cor} \label{algCVcor}
For every $p \in (1,\infty)$, there exists a positive constant $C'(p)$ such that for $t \geq 0$, we have
\begin{equation} \\
\Vert S_t^{*}\mu_1-S_t^{*}\mu_2 \Vert^*_L \leq C't^{-1/13p},\qquad t \geq 1,
\end{equation}
for any probability measures $\mu_1$, $\mu_2$ on $L_p$.
\end{cor}

\indent
Note that all the estimates in the previous sections still hold for a stationary solution, since they hold uniformly for any initial condition in $L_1$ for large times and a stationary solution has time-independent statistical properties. It follows that those estimates still hold when averaging in time and in ensemble (denoted by $\lbrace \cdot \rbrace$) is replaced by averaging solely in ensemble, i.e. by integrating  with respect to $\mu$. Namely, Theorem~\ref{whiteavoir}, Theorem~\ref{whiteavoir2} and Theorem~\ref{whitespectrinert} imply, respectively, the following results.

\begin{theo}
For $m \in \lbrace 0,1 \rbrace$ and $p \in [1,\infty]$, or for $m \geq 2$ and $p \in (1,\infty]$,
\begin{equation} \nonumber
\Big( \int {\left|u\right|_{m,p}^{\alpha} d \mu(u)} \Big)^{1/\alpha} \overset{m,p,\alpha}{\sim} \nu^{-\gamma},\quad \alpha>0.
\end{equation}
\end{theo}

\begin{theo}
For $\alpha \geq 0$ and $\ell \in J_1$,
$$
\int{ \Big( \int_{S^1}{|u(x+\ell)-u(x)|^p dx} \Big)^{\alpha} d \mu(u)} \overset{p,\alpha}{\sim} \left\lbrace \begin{aligned} & \ell^{\alpha p},\ 0 \leq p \leq 1. \\ & \ell^{\alpha p} \nu^{-\alpha (p-1)},\ p \geq 1. \end{aligned} \right.
$$
On the other hand, for $\alpha \geq 0$ and $\ell \in J_2$,
$$
\int{ \Big( \int_{S^1}{|u(x+\ell)-u(x)|^p dx} \Big)^{\alpha} d \mu(u)} \overset{p,\alpha}{\sim} \left\lbrace \begin{aligned} & \ell^{\alpha p},\ 0 \leq p \leq 1. \\ & \ell^{\alpha},\ p \geq 1. \end{aligned} \right.
$$
\end{theo}

\begin{theo}
For $k$ such that $k^{-1} \in J_2$, we have:
$$
\int{ \frac{\sum_{|n| \in [M^{-1}k,Mk]}{|\hat{u}(n)|^2}}{\sum_{|n| \in [M^{-1}k,Mk]}{1}} d \mu(u)} \sim k^{-2}.
$$
\end{theo}

%

\section*{Acknowledgements}
\indent
I am very grateful to A.Biryuk, U.Frisch, K.Khanin, S.Kuksin and A.Shirikyan for helpful discussions. A part of the present work was done during my stay at the AGM, University of Cergy-Pontoise, supported by the grant ERC 291214 BLOWDISOL: I would like to thank all the faculty and staff, and especially the principal investigator
\\
F.Merle, for their hospitality.
\bigskip
\begin{center}
Alexandre Boritchev
\\
Laboratoire AGM
\\
University of Cergy-Pontoise
\\
2 av. Adolphe Chauvin 
\\
95302 CERGY-PONTOISE CEDEX
\\
FRANCE
\end{center}

\bibliographystyle{plain}
\bibliography{Bibliogeneral}

\end{document}